\documentclass[a4paper,11pt]{article}

\usepackage{graphicx}
\usepackage[all]{xy}
\usepackage[english]{babel}

\usepackage{amsmath,amssymb}

\usepackage{amsthm}

\usepackage{amscd}

\usepackage{mathrsfs}

\begin{document}

\title{Combinatorial aspects of nodal curves}
\author{Simone Busonero, Margarida Melo, Lidia Stoppino\footnote{The authors thank the University of Catania and, in particular, the organizers of PRAGMATIC 2004, for offering them the possibility and the conditions for developing this joint work.}}

\newcommand{\av}{``}
\newcommand{\Z}{\mathbb Z}
\newcommand{\forn}{\forall n\in\mathbb{N}}
\newcommand{\raw}{\to}
\newcommand{\Raw}{\Rightarrow}
\newcommand{\law}{\gets}
\newcommand{\Law}{\Leftarrow}
\newcommand{\Lra}{\Leftrightarrow}

\newtheorem{teo}{Theorem}[section]
\newtheorem{prop}[teo]{Proposition}
\newtheorem{lem}[teo]{Lemma}
\newtheorem{rem}[teo]{Remark}
\newtheorem{defi}[teo]{Definition}
\newtheorem{ex}[teo]{Example}
\newtheorem{cor}[teo]{Corollary}

\maketitle
\begin{abstract}
To any nodal curve $C$ is associated the degree class group, a combinatorial invariant which plays an important role in the compactification of the generalised Jacobian of $C$ and in the construction of the N\'eron model of the Picard variety of families of curves having $C$ as special fibre. In this paper we study this invariant. More precisely, we construct a wide family of graphs having cyclic degree class group and we provide a recursive formula for the cardinality of the degree class group of the members of this family. Moreover, we analyze the behaviour of the degree class group under standard geometrical operations on the curve, such as the blow up and the normalisation of a node.
\end{abstract}

\tableofcontents

\section*{Introduction}

Let $C$ be a reduced nodal curve defined over an algebraically closed field $k$.
Let $f:\mathcal C\rightarrow B =\mbox{Spec }R$, where $R$ is a discrete valuation ring with residue field $k$, be a family of generically smooth 
nodal curves, such that  the special fibre is isomorphic to $C$. 
Consider the set of Cartier divisors $D$ on $\mathcal C$ supported on $C$; the associated line bundles $\mathcal O_\mathcal C(D)$ are called {\em twisters}.
Given a line bundle $\mathcal L$ on $\mathcal C$, the line bundles of the form
$\mathcal L\otimes \mathcal O_\mathcal C(D)$ clearly agree with $\mathcal L$ on the
general fibre, but differ on the special one (also the multi-degrees on $C$ are different).
Conversely, if a line bundle $\mathcal M$ agrees with $\mathcal L$ on the general fibre, 
then it has to be of the form $\mathcal M=\mathcal L\otimes \mathcal O_\mathcal C(D)$.
So, if we consider the Picard functor $\mathcal{P}ic_f$ of the family $f$, we can say that {\em the twisters cause the nonseparatedness} of this functor. 
If the total space $\mathcal C$ is regular, the multidegrees of the twisters depend only on the combinatorics of $C$ (if $\mathcal C$ is only normal, one has to consider the type of rational singularities it has).

The object of study of this paper is the group of classes of multidegrees on $C$ modulo the multidegree of twisters, the so-called {\em degree class group} of $C$, DCG for short (see section \ref{pre} for the precise definition).
It is clearly a purely combinatorial invariant of the curve.
In this form it was introduced in \cite{cap}, in order to describe and handle the fibres of the compactification of the universal Picard variety $\overline P_{d,g}$ over the moduli space of stable curves $\overline M_g$. 
In particular, the fibre of $\overline P_{d,g}$ over $[C]\in \overline M_g$ can be seen as a compactification of the generalised Jacobian $J_C$, and there is an injective map between the set of its irreducible components and the DCG of $C$.

\medskip

In fact, the DCG associated to a nodal curve has been extensively studied in Arithmetic Geometry as a particular case of the following more general construction (see for instance  \cite{ray}, \cite{BLR} and \cite{lorjac}). 
Given a discrete valuation ring $R$ with residue field $k$ (not necessarily algebraically closed)  and quotient field $K$,  let $X\rightarrow \mbox{Spec }R$ be a flat projective curve such that $X$ is regular and the generic fibre $X_K$ is geometrically irreducible. Then, under some technical assumptions, there exists a N\'eron model for the jacobian $J_K$ of of $X_K$. 
The special fibre $X_k$ is of the form $\sum m_i X_i$ with $X_i$ irreducible and distinct; there is a natural group  $ \Phi$ associated to the intersection matrix of the $X_i$'s, which is the \emph{ group of connected components of the special fibre of the N\'eron model}. 
When $X_k$ is a nodal curve (in particular $m_i=0$ for every $i$), $\Phi$ coincides with the DCG of $X_k$. 
The structure of the  group $\Phi$ has been the object of a series of papers by Lorenzini (\cite{lorgraphs}, \cite{lorjac}, \cite{Lorarit}, \cite{Lorlap}, \cite{Lorgroup}).

More recently,  Caporaso in  \cite{capneron} gave a geometric counterpart of this construction, showing the existence of a space over $\overline M_g$ such that for every regular family $f:\mathcal X\rightarrow B$ of stable curves the N\'eron model of the Picard variety of degree $d$ of $\mathcal X$ is obtained by base change via the moduli map $B\rightarrow \overline M_g$.

\medskip

Another incarnation of the  DCG is in Combinatorics, as an invariant of graphs (see for instance \cite{biggs}, \cite{biggschip}, \cite{tatiana}) and in this field it goes under many other names, such as critical group, determinant group, Picard group, Jacobian group.  
Also from the point of view of Combinatorics, a typical problem is to compute the structure of this group. 
It has been solved completely only for a few families of graphs.
The family of graphs with cyclic DCG constructed in section \ref{ciclici} is a new contribution in this sense. 

\medskip

It is clear from the  above exposition that the DCG of a nodal curve comes out as a significant invariant of the curve  in many geometric contexts. It is therefore natural to ask if it is possible to \emph{classify} nodal curves using their DCG. 
In particular, one could hope to use this discrete invariant to try and \emph{stratify} the moduli space of stable curves $\overline M_g$. 
As the DCG is in fact an invariant of the dual graph of the curve, it could give a coarser stratification than the one given by topological type. 
Moreover, this stratification would be extremely different from the one associated to the number of nodes (just observe that any compact type and any irreducible curve, regardless of the number of nodes, have trivial DCG). 

This was in fact the original motivation of this work. 
However, as the numerous results both in Arithmetics and in Combinatorics clearly show, this task is far too ambitious; for instance, also the problem of classifying all curves having cyclic DCG is very hard to solve. 
Keeping in mind the idea of a classification associated to the DCG, in this paper we try and improve the understanding of the connections between the geometrical properties of a nodal curve and the structure of its DCG.
We perform the computation of the order and structure of the DCG for some types of curves, and compute several examples.
Moreover, we study the relation between geometrical operations on the curve and the corresponding modifications on the DCG, giving  some useful formulas. 

More precisely, the contents of the paper are the following:
In the first section we introduce the main objects and techniques of our study, and we present a proof of the equality between the cardinality of the DCG of a curve and the complexity of its associated dual graph (Kirkoff's Matrix Tree Theorem).

In section \ref{sec2}, after studying some simple cases, we describe a family of graphs having  cyclic degree class group. 
Moreover, we list all the possible graphs for \emph{stable} curves of genus $2$ and $3$.
 
In section \ref{sec3} we analyse the behaviour of the DCG under the geometric operations of blow up, normalisation and smoothing of a node. This is a problem arising in the geometric applications of the DCG. 
One example is the following: if we consider  families of nodal curves (again with general smooth fibre) such that the total space is {\em normal}, the nodes $P_1,\cdots P_n$ of the special fibre $C$ will correspond to rational singularities of the total space, say of type $A_{m_1}, \cdots A_{m_n}$. The group of components of the N\'eron model of the relative Jacobian is not the DCG of $C$, but the DCG of the blow up of $C$ $m_i$ times in the $i$-th node, respectively.

In section \ref{blow} we translate a standard graph theory result in terms of geometric operations on the curve. 
This way we can obtain a general formula (Theorem \ref{perepe}) relating the DCG of a blown up curve to the DCG of its partial normalisations.

In section \ref{vine} we obtain some results on the cardinality and structure of the so-called vine curve, i.e. a curve with two smooth components meeting in $N$ nodes.
Our computations, although obtained with different techniques, can be derived from the results of  \cite{BLR} (proposition 9.6.10) and from \cite{Lorarit} (example 2.5 and successive claims). 

\bigskip

\noindent{\bf Acknowledgements} We wish to express our deep gratitude to Lucia Caporaso for the suggestion of the problem that led us to this work, that she patiently  supervised, and to Cinzia Casagrande for precious remarks and corrections on the preliminary version of the paper. The second author thanks Jo\~ao Gouveia for improving some calculations and the third author thanks Ludovico Pernazza for his patient help with some unpleasant computations. 
Moreover, we wish to thank Dino Lorenzini for having pointed out several inaccuracies in the previous version of the paper, and for the kind interest he showed towards our work.


\section{Preliminaries and first results}\label{pre}

Let $k$ be an algebraically closed field.
Throughout the  paper a curve will mean a connected reduced nodal curve projective over $k$.
The genus $g=g(C)$ of a nodal curve $C$ is the arithmetic genus $h^0(\omega_C)$, where $\omega_C$ is the dualising sheaf of $C$.
For each such curve $C$ we will call $\gamma(C)$ the number of irreducible components of $C$ and 
$\delta(C)$ the number of nodes of $C$.

\subsubsection*{The dual graph of a curve}

To a curve $C$ we can associate a graph $\Gamma_C$, i.e. a symplicial complex of dimension at most 1, called the \textit{dual graph}, in the following way:
\begin{itemize}
\item to each irreducible component $A$ corresponds a vertex $v_A$ (i.e. a $0$-dimensional symplex); 
\item to each node intersecting the components $A$ and $B$ (where $A$ and $B$ can coincide) 
corresponds an edge (1-dimensional symplex) connecting the vertices $v_A$ and $v_B$.
\end{itemize}
Thus $\Gamma_C$ has $\gamma(C)$ vertices (i.e. it has {\em order} $\gamma(C)$), $\delta(C)$ edges, 
and among the edges there is a loop 
for every node lying on a single irreducible component of $C$. 
Note that two vertex can be joined by more than one edge.

Recall that the first Betti number of $\Gamma_C$ is 
$$b_1(\Gamma_C):= \delta (C)-\gamma (C)+1.$$
(in the general formula, 1 is substituted by the number of connected components of $C$).

Recall that, for any nodal curve $C$ if  $C_1,\ldots,C_\gamma$ are its irreducible components, and   $g_i=g(C_i)$, then the arithmetic genus of $C$ is
$$g=\sum_{i=1}^{\gamma}g_i+\delta(C)-\gamma(C)+c,$$
where $c$ is the number of connected components of $C$ and $\delta$ is the number of nodes of $C$.
Notice that, as we consider all curves to be connected, in what follows we will always use $c=1$.

We can also construct a {\em weighted} graph, associating to any vertex $v$ the
genus $g_v$ of the corresponding component. 
In fact the weighted graph constructed this way encode all the topological 
information about the curve. 

\begin{rem}
\upshape{
Observe that {\em every} connected graph can be considered as the dual graph of a curve.}
\end{rem}

\subsubsection*{Complexity of a graph}

\begin{defi}
Let $\Gamma$ be a graph. A {\it spanning tree} of $\Gamma$ is a subgraph of $\Gamma$ which is a tree having 
the same vertices as $\Gamma$. The complexity of $\Gamma$, indicated by the symbol $c(\Gamma)$, is the number of spanning trees contained 
in $\Gamma$.
\end{defi}
Not every introductory book on graph theory treats this topic. See for reference \cite{biggs}, Section 6, \cite{berge}, cap.3 $\natural$ 5 and \cite{west}, Section 2.2.

Observe that $c(\Gamma)=0$ if and only if $\Gamma$ is not connected, and that if $\Gamma$ is a connected tree  $c(\Gamma)=1$.

For the complexity of the dual graph associated to a curve $C$, we will often use the symbol $c(C)$, instead of
$c(\Gamma_C)$.

\subsubsection*{Degree class group}

Let $\left\{C_i\right\}_{i=1,...,\gamma}$ be the irreducible components of a curve $C$. Define 
$$
\begin{matrix}k_{ij}:= & \left\{\begin{array}{l}\;\;\, \sharp(C_i\cap C_j) \, \mbox{ if }i\not= j\\  \\
-\sharp(C_i\cap \overline{C\setminus C_i}) \, \mbox{ if } i=j\\\end{array}\right.\end{matrix}
$$
As $C_i\cap \overline{C\setminus C_i}=\bigcup_{j\not= i} C_i\cap C_j$, we have that for fixed $i$, $\sum_j k_{ij}=0$.
For every $i$ set 
$$\underline{c}_i:=(k_{i1},\ldots,k_{i\gamma }) \in \mathbb Z^\gamma.$$ 
Call $Z:=\{\underline{z}\in\mathbb{Z}^\gamma : |\underline z|=0\}$.
As observed before, $\underline c_i\in Z$. Let us call $\Lambda_C$ the sublattice of $Z$ spanned by 
$\{\underline c_1,\ldots,\underline c_\gamma\}$. In fact, $\Lambda_C$ is a lattice in  $Z$ (it
has rank $\gamma -1$) as we will show in a moment (see \cite{capneron} for a geometric proof of this fact).

\begin{rem}\label{generatori}
\upshape{Fix a one-to-one correspondence between the set $V$ of vertices of the graph and the elements of the 
canonical basis of $\mathbb Z^\gamma$, and call $e_v$ the element of the basis associated to $v$ 
with respect to the correspondence chosen; observe that, for any $w\in V$, $Z$ is generated by the elements 
$\{e_w-e_v, \,v \in V\}$.}
\end{rem}

\begin{defi}
The degree class group of $C$ is the finite abelian group $\Delta_C:=Z/\Lambda_C$.
\end{defi}
For short, we will denote the degree class group as DCG.
This name was given in \cite{cap} where such a group was introduced to compactify the generalised 
Jacobian of stable curves.

\begin{rem}\upshape{
It is important to notice that the DCG depends {\em only} on the dual graph
of the curve: clearly we can define it  for any  graph. Indeed, given a loopless connected graph $\Gamma$ with vertices
$\{ v_1,...,v_\gamma\}$, we simply define the $k_{ij}$'s in the following way:

$$
\begin{matrix}k_{ij}:= & \left\{\begin{array}{l}\;\;\, \sharp\{\mbox{edges connecting }v_i\mbox{ and }v_j \}\, \,
\mbox{ if }i\not= j\\  \\
- \sharp \{\mbox{edges touching }v_i\}=- (\mbox{degree of }v_i)\,\, \mbox{ if } i=j\\\end{array}\right.\end{matrix}
$$
We will call $\Delta_\Gamma$ the DCG associated to the graph $\Gamma$. For general connected graphs, we define the DCG 
as the DCG of the corresponding loopless graph.}
\end{rem}

Let $M$ be the $\gamma\times\gamma$ matrix whose columns are the $\underline c_i$'s. We will call $M$ the {\em intersection matrix} \footnote{Readers familiar with graph theory can observe that $M$ is obtained 
from the adjacency matrix subtracting the vertex degrees on the diagonal. This matrix is frequently referred as the Laplacian of the graph (see for instance \cite{Lorlap})}.

The following theorem, known as Kirkoff's Matrix Tree Theorem, will be a key ingredient for our analysis of the DCG. Given its importance, we present here 
also a proof. See for reference \cite{west}. There are at least other two proofs of this theorem:
see \cite{godsil-royle} and  \cite{chaiken}.

\begin{teo}{\upshape (Matrix Tree Theorem)}
Let $s, t \in \{1,\ldots \gamma\}$.
Using the above notations, if $M^s_t$ is obtained by $M$ by deleting the $t$-th column and the $s$-th row, then 
$$c(\Gamma)=(-1)^{s+t+\gamma -1}\mbox{det}(M^s_t).$$
\end{teo}
\begin{proof}
The sum of the columns of $M$ is zero, thus when we replace the $s$-th column of $M^s_t$ with the $t$-th column of the 
matrix obtained from $M$ by deleting the $s$-th row, the sign of the determinant of $M^s_t$ is reversed, 
whereas its absolute value remains unchanged. 
Successively, we can permutate the columns so that the matrix becomes the one obtained by $M$ by deleting the s-th column and the s-th row. 
The sign of this permutation is $(-1)^{\mid s-t\mid -1}$.
Therefore,  
$$\det M^s_s=(-1)^{s-t}\det M^s_t$$
so we can suppose $s=t$. Then we have to prove that 
\begin{equation}\label{thesis}
\det (-M^t_t)=c(\Gamma)
\end{equation}
holds for every $t=1$,\dots $\gamma-1$.

From now on, fix an orientation on the graph and an enumeration on its edges. 
Let $I$ be the incidence matrix of $\Gamma$: 
the entries of $I$ are $a_{i,j}=1$ when $v_{i}$ is the tail of $e_{j}$, $a_{i,j}=-1$ when $v_{i}$ is the head of $e_{j}$ and $a_{i,j}=0$ otherwise.
Observe that $-M=I\cdot I^{T}$.

Let $I^{\star }$ be the result of deleting row $t$ of $I$, so $-M^t_t =I^\star\cdot (I^{\star })^{T}$.
The Binet-Cauchy formula computes the determinant of a product of non-square matrices using the determinants of maximum square submatrices of the factors:
let $A$ be $p\times m$, let $B$ be $m\times p$, $m\geq p$, then 
$det(AB)=\sum _{\mid S\mid =p}A_{S}B_{S}$, 
where $A_{S}$ is the submatrix of $A$ consisting of the columns indexed by $S$ and $B_{S}$ is the submatrix of $B$ 
consisting of the rows indexed by $S$. 
Since $I^{\star }$ is $(\gamma -1)\times \delta$ and $\Gamma$ is connected (and so $\gamma -1=\delta -g\leq \delta$), we can apply the Binet-Cauchy formula to $-M^t_t=I^{\star }\cdot (I^{\star })^{T}$, so that $S$ runs over all the sets of $\gamma -1$ edges of $\Gamma$, $A_{S}$ is a $(\gamma -1)\times (\gamma -1)$ submatrix of $I$ and $B_{S}$ is $A_{S}^{T}$, so
$$\mbox{det}(-M^t_t)=\sum _{S} (\mbox{det}I^{\star } _{S})^{2},$$
where the sum runs over all the sets of $\gamma -1$ edges of $\Gamma$.  

We will prove below that the determinant of every $(\gamma -1)\times (\gamma -1)$ submatrix of $I$ is $\pm 1$ if the associated  set of $\gamma -1$ edges form a spanning tree of $\Gamma$ (point 1), while it is zero otherwise (point 2). 
Observe that if we assume this, the absolute value of the previous summand counts exactly all the possible spanning trees in $\Gamma$ and so we obtain formula (\ref{thesis}).
 
\smallskip 
 
1) In the first case we use induction on $\gamma$.
For $\gamma =1$, it's clear because by convention a $0\times 0$ matrix has determinant 1.
For $\gamma >1$, let $T$ be a spanning tree whose edges are columns of a $(\gamma -1)\times (\gamma -1)$ submatrix $B$ of 
$I$. 
Since the sum of the degrees of the vertices is two times the number of the edges, a tree has at least two leaves, i.e. 
vertices whose degree is 1, and since only one row of $I$ is deleted, $B$ has a row corresponding to a leaf $v$ of $T$. 
This row has only one nonzero entry in $B$, which is $\pm 1$;  when one computes the determinant by expanding along this 
row, the only submatrix $B^{\prime }$ with nonzero coefficient corresponds to the spanning subtree of $\Gamma -v$ (obtained by deleting $v$ and its incident edge from $T$). We can therefore apply the inductive hypothesis to $B^{\prime }$.

\smallskip

2) Now, suppose that the $\gamma -1$ edges corresponding to the columns of $B$ do not form a spanning tree. 
Then they contain a cycle $C$. Indeed, if this were not the case, calling $\Gamma'$ the subgraph made of this edges, we would have
$$0=b_1(\Gamma' )=E(\Gamma ')-V(\Gamma ')+\sharp \pi _{0}(\Gamma ')\geq V(\Gamma) -V(\Gamma ')\geq 0,$$ 
where $\pi _{0}$ is the set of connected components, $E$ and $V$ the number of edges and of vertices respectively. 
Therefore $V(\Gamma ') =V(\Gamma )=E(\Gamma ')+1$ and $\sharp \pi _{0}(\Gamma ')=1$, so that the edges would form a spanning tree.

We form a linear combination of the columns in this way: with coefficient $0$ if the corresponding edge is not in $C$, $+1$ if it is 
followed forward by $C$, and $-1$ if it is followed backward by $C$.
The result is of total weight $0$ at each vertex, so the columns are linearly dependent, which yields det$B=0$.
\end{proof}

The Matrix Tree Theorem assures that $M$ has rank $\gamma-1$, i.e. that $\Lambda_C$ is indeed a lattice.
Moreover, it allows us to relate the cardinality of the DCG of a 
curve $C$ with the complexity of its dual graph, as we see below.

For $r\in\{1,\ldots,\gamma\}$, consider the isomorphism $\alpha_r:Z\stackrel{\sim}\longrightarrow\mathbb{Z}^{\gamma - 1}$
which consists of deleting the $r$-th component.
The group $\Delta_C$ is the quotient of $\mathbb{Z}^{\gamma - 1}$ by the lattice generated by 
$$\underline c_i^\prime:=(k_{1i}
,\ldots,\widehat{k_{ri}},\ldots,k_{\gamma i}).$$
Observe that again $\sum_i \underline c_i^\prime = \underline 0 \in \mathbb{Z}^{\gamma -1}$.
Therefore $\Delta_C$ is presented by the matrix $M^\star$ obtained from $M$ deleting a column and the $r$-th row
(for presentation of modules by integer matrices see \cite{ar}).
Consider now the following sequence
$$\mathbb{Z}^{\gamma-1}\stackrel{M^\star}\longrightarrow \mathbb Z^{\gamma -1}\longrightarrow \Delta_C\longrightarrow 0$$
where the first map is the linear map associated to $M^\star$.
By diagonalisation of integer matrices (cf. \cite{ar}), there exists a diagonal presentation matrix $D$ for $\Delta_C$, 
i.e. there exists $P,Q\in GL(\gamma - 1,\mathbb Z)$ and a diagonal matrix $D\in \mbox{Mat}(\gamma - 1,\mathbb Z)$ 
such that
$$P M^\star Q^{-1}=D.$$
The absolute values of the entries on the diagonal of $D$ correspond to the order of the cyclic factors of $\Delta_C$ (the so-called invariant factors; notice that this is in fact the structure theorem for abelian groups). \footnote{
Notice that although these matrices are diagonalisable also in $\mathbb R$ (being symmetric), the eigenvalues do not correspond at all to the invariant factors, not even in the case they are integer; a nice counterexample can be found in Section 9.2 of \cite{tatiana}.}
Therefore,
$$\sharp \Delta_C=|\mbox{det}(D)|=|\mbox{det}(P)\mbox{det}(M^\star)\mbox{det}(Q^{-1})|=|\mbox{det}(M^\star)|.$$
So we can conclude that {\em the cardinality of the DCG of a curve $C$ is the complexity of the dual graph $\Gamma_C$.}

It's worth noticing that this equality is well known; see for instance \cite{Lorgroup}, remark on pag. 280.
In \cite{OS} and in \cite{cap} there is a proof involving a \av cohomological" computation of the DCG and a theorem of Kirkoff-Trent. 


\section{Computing the cardinality and  the structure of the DCG}\label{sec2}

We have seen in the previous section that given a curve, we can find the cardinality of its DCG 
simply by computing a determinant, and the structure of
its DCG performing a diagonalization of integer matrices. 
A natural question arising at this point is the following: 
what kind of curves have fixed DCG, or DCG with some fixed properties, i.e. can we somehow {\em classify} curves 
using this invariant?
The results contained in this section, or even in the whole paper, can be seen as evidences of the fact that this is a very complicated and involved problem.

In this section we compute several examples, and we state some partial results about curves whose DCG is cyclic. 

\bigskip

Let us start by considering the simplest situations.
For example, what kind of curves have DCG trivial? Clearly this
means that the dual graph associated to $C$ is a tree, once removed all the possible loops it may 
have.
Therefore $C$ must be such that any non disconnecting node has both preimages in the same 
component of the normalisation. 

\begin{rem}
{\upshape 
Clearly, to remove or to attach to one vertex of a graph another graph with complexity 1 
doesn't change the complexity.
On the other hand, notice that it does change the associated curve. 
From now on in this section, we will consider 
graphs modulo this operation.}
\end{rem}
\noindent Here we list the possible loopless graphs (modulo trees) with complexity $2$, $3$, $4$:

\bigskip

\begin{center}

\begin{tabular}{|c|c c c c|}

\hline
complexity 2 
&
 \xymatrix@R=.2pc{
*{\bullet} \ar @{-} @/_/[r]  & *{\bullet} \ar@{-} @/_/[l]\\
{}
}
&
\hspace{.9cm}
&

&

\\
\hline

complexity 3 
&
\xymatrix@R=.5pc{
& \\
*{\bullet} \ar @{-} @/_/[r] \ar@{-}[r] & *{\bullet} \ar@{-}@/_/[l]\\
{}\\
{}
}
&
\xymatrix@=.6cm{
& *{\bullet} \ar@{-}[dl] \ar@{-}[dr]& \\
*{\bullet} \ar@{-}[rr] & & *{\bullet}
}
&

&

\\

\hline

complexity 4
&
\xymatrix@R=.7pc{
\\
*{\bullet}  \ar @{-} @/_.3pc/[r] \ar @{-}@/_.8pc/[r]/& *{\bullet} \ar@{-}@/_.3pc/[l] \ar@{-}@/_.8pc/[l]
}
&
\xymatrix@R=.5cm{\\
*{\bullet} \ar @{-} @/_/[r]  & *{\bullet} \ar@{-} @/_/[l] \ar @{-} @/_/[r]  & *{\bullet} \ar@{-} @/_/[l]\\
{}
}
&
\xymatrix@=.6cm{
& *{\bullet} \ar@{-}[dl] \ar@{-}[dr]& \\
*{\bullet} \ar@{-}[rr] \ar@{-}@/_/[rr]& & *{\bullet}
}
&
\xymatrix{
*{\bullet} \ar@{-}[r] \ar@{-}[d] & *{\bullet} \ar@{-}[d] \\
*{\bullet} \ar@{-}[r] & *{\bullet}
}\\
\hline
\end{tabular}

\end{center}

\bigskip

\begin{rem}\label{prodotto}{\upshape (cf. also \cite{Lorgroup})
If  $\Gamma$ be a graph obtained attaching  graphs $\Gamma_1$ and $\Gamma_2$ in one vertex.
Then $$\Delta_\Gamma=\Delta_{\Gamma_1}\oplus\Delta_{\Gamma_2}.$$
Indeed, let $n$ be the order of $\Gamma$, $k$ the order of $\Gamma_1$. 
Choose an ordering of the vertices of $\Gamma$ such that the first $k$ belong to $\Gamma_1$ (so the vertex
of index $k$ is the common vertex of $\Gamma_1$ and $\Gamma_2$). Let $M$ be the intersection matrix of $\Gamma$
with respect to this ordering. Observe that if we remove the $k$-th row and column 
from $M$ we obtain a block matrix, and apply the Matrix Tree Theorem.}
\end{rem}

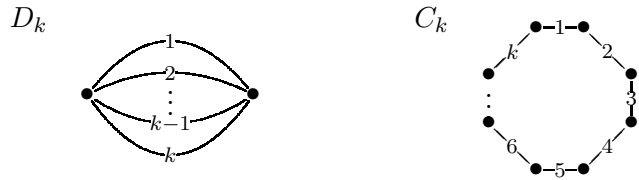
\begin{figure}[!htp]
\begin{center}
\begin{tabular}{rcccc}
$D_k$
&
\xymatrix@R=1pc{\\
*{\bullet}  \ar @{-} @/_.9pc/[rr]|-{k-1} \ar @{-}@/_1.9pc/[rr]|-{k}/& \vdots &*{\bullet} \ar@{-}@/_.7pc/[ll]|-{2} \ar@{-}@/_1.7pc/[ll]|-{1}
}

&
\hspace{1.5cm}

$C_k$
&
\xymatrix@=1pc{
& *{\bullet} \ar@{-}[r]|-{1} & *{\bullet} \ar@{-}[dr]|-{2} & \\
*{\bullet} \ar@{-}[ur]|-{k} & & & *{\bullet} \ar@{-}[d]|-{3} \\
*{\bullet} \ar@{}[u]|-{\vdots} \ar@{-}[dr]|-{6} & & & *{\bullet} \ar@{-}[dl]|-{4} \\
& *{\bullet} \ar@{-}[r]|-{5} & *{\bullet} &
}
\end{tabular}
\caption{Two vertices attached by $k$ edges, $D_k$, and the $k$-cycle, $\mathcal C_k$.}\label{figura}
\end{center}
\end{figure}

\begin{ex}\label{dk}
{\upshape Call $D_k$ the graph made of two vertices attached by $k$ edges showed in fig. \ref{figura} (this is the graph of a vine curve, as defined in section \ref{vine}).
The intersection matrix is $\begin{pmatrix}-k & k\\
k &-k\\ \end{pmatrix}$, so clearly $\Delta_{D_k}\cong \mathbb Z/k\mathbb Z$.}
\end{ex}

\begin{ex}
{\upshape Call $\mathcal C_k$ the $k$-cycle (fig. \ref{figura}). Using the definition of complexity it is easy to 
see that the cardinality of its DCG is $k$. Ordering clockwise the vertices, we have
$$
\underline c_i=e_{i-1}-2e_i+e_{i+1},
$$
where the indexes are obviously considered mod $k$.
Therefore 
$$e_i-e_{i+1}=\underline c_{i+1}+ e_{i+1}-e_{i+2},$$ 
so $\Delta_{\mathcal C_k}$ has one generator (remember Remark \ref{generatori})
and again we can conclude that the DCG is isomorphic to $\mathbb Z/k\mathbb Z$.
}
\end{ex}

One of ours key tools is the following well-known result.

\begin{prop}\label{comb}
Let $\Gamma$ be a graph. If $e$ is an edge of $\Gamma$ which is not a loop,
call $\Gamma -e$ the graph obtained from $\Gamma$ removing $e$, and $\Gamma\cdot e$ the one obtained  contracting $e$. 
Between the complexities of these three graphs the following relation holds:
\begin{equation}\label{fond}
c(\Gamma)=c(\Gamma - e)+c(\Gamma\cdot e).
\end{equation}
\end{prop}
\begin{proof}
Just observe that the spanning trees of $\Gamma\cdot e$ correspond bijectively to the spanning trees of $\gamma$ containing $e$, while the spanning trees of $\Gamma - e$ are clearly the ones of $\Gamma$ not containing $e$.
\end{proof}
Let $v$ and $w$ be two vertices of $\Gamma$  having exactly $r$ edges $\{e_1,\ldots ,e_r\}$ in common. 
Let $\Gamma'$ be the graph obtained from $\Gamma -\{e_1,\ldots ,e_r\}$ by identifying $v$ and $w$.
From the above proposition, it follows easily by induction the formula 
$$c(\Gamma)= c(\Gamma -\{e_1,\ldots ,e_r\}) + r c(\Gamma '),$$
which is established, and extensively exploited, also in  \cite{Lorarit}.

\subsection{A family of graphs with cyclic DCG}
A natural question to ask is whether is possible to classify all graphs whose DCG is cyclic. 
Even if they seem to be very different, we see below that the two examples above are particular cases of a more general type of graphs.

\bigskip

Let $n$ be a positive integer.
Let $\underline{k}$ be an element of $(\Z _{>1})^{n}$ and $\underline{h}$ be an element of $(\Z_{>0}) ^{n}$ such that the $i$-th coordinate of $\underline{h}$ is smaller than the 
$i$-th coordinate of $\underline{k}$. 
For each coordinate $k_{j}$ of $\underline{k}$, we assign a $k_{j}$-cycle $\mathcal{C} _{k_{j}}$ whose set of vertices is a 
double indexed set $\{ v^{j}_{1},\dots ,v^{j}_{k_{j}} \}$ ordered clockwise.
Then to a coordinate $h_{j}$ of $\underline{h}$ corresponds a vertex $v^{j}_{h_{j}}\in \mathcal{C} _{k_{j}}$.
Given the data $n$, $\underline{k}$, $\underline{h}$, we will build a graph $CS^{n}(\underline{k};\underline{h})$, 
using induction on $n$. 

For $n=1$, $\underline{k}=k$, $\underline{h}=h$, we define $CS^{1}(k;h)\colon = \mathcal{C} _{k}$ (this way we obtain all the cycles).

For $n=2$, define a set-map $A_{2}$ from a subset of $V(\mathcal{C} _{k_{2}})$ to a subset of $V(CS^{1}(k_{1};h_{1}))$, by $A_{2}(v^{2}_{1})=v^{1}_{h_{1}}$, $A_{2}(v^{2}_{k_{2}})=v^{1}_{h_{1}+1}$.
Then 
$$V(CS^{2}(\underline{k};\underline{h}))\colon = V(CS^{1}(k_{1};h_{1}))\sqcup _{A_{2}}V(\mathcal{C} _{k_{2}})$$
$$E(CS^{2}(\underline{k};\underline{h}))\colon = E(CS^{1}(k_{1};h_{1}))\cup E(\mathcal{C} _{k_{2}})/\{ v^{1}_{h_{1}}v^{1}_{h_{1}+1}\} \cup \{ v^{2}_{k_{2}}v^{2}_{1}\}$$

The proof of the inductive step is analogous to step $n=2$.

\medskip

\begin{figure}[!htp]
\begin{center}
\xymatrix@C=3pc{
 &&&&& *{\bullet} \ar@{-}[dl]_>{v_2^3=v_1^4} \ar@{-}[dr]^<{v_2^4} &&&&\\
 & *{\bullet} \ar@{-}[r]^<{v_2^1=v_1^2} \ar@{-}[dd] & *{\bullet} \ar@{-}[rr]^<{v_2^2=v_1^3} \ar@{-}[dd] & & *{\bullet} \ar@{-}[dd] & & *{\bullet} \ar@{-}[dd] \ar@{-}@/_/[rrdd] \ar@{-}@/^1pc/[rrdd]^<{\;\;\;\:\: v_3^4=v_1^5=v_1^6} & &  \\
*{\bullet} \ar@{-}[ur]^<{v_1^1} \ar@{-}[dr]_>{v_3^1=v_4^2} &&&&& && & \\
 & *{\bullet} \ar@{-}[r]_>{\hspace{-.2cm}v_3^2=v_5^3}  & *{\bullet} \ar@{-}[dr]_>{v_4^3}& & *{\bullet}  & & *{\bullet} \ar@{-}[rr]_>{v_2^5=v_2^6}&& *{\bullet}\\
& & & *{\bullet} \ar@{-}[ur]^>{\vspace{.2cm}v_3^3=v_5^4} && *{\bullet} \ar@{-}[ur]_<{v_4^4}_>{v_3^4=v_3^5} \ar@{-}[ul]
}
\caption{A representation of $CS^6((3,4,5,6,3,2);(2,2,2,3,1,1))$.}\label{polygonal}
\end{center}
\end{figure}
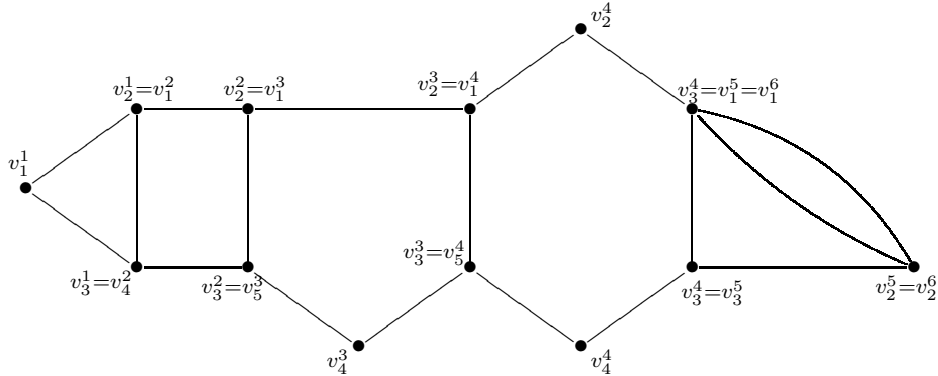

We can draw $CS^{n}(\underline{k};\underline{h})$ as a chain of polygons such that each polygon and the following one 
are attached at only one edge (see figure \ref{polygonal} for an example).
Therefore the graph $D_e$ of Example \ref{dk} is isomorphic to $CS^{e-1}(\underline{2};\underline{1})$.
 
\begin{teo}\label{ciclici}
The degree class group of $CS^{n}(\underline{k};\underline{h})$ is cyclic.
\end{teo}
\begin{proof} Let us fix an ordering of the $\gamma$ vertices of $CS^n(\underline{k}, \underline{h})$. So, for each vertex $v_i^j$ of $CS^n(\underline{k}, \underline{h})$, let $e_{v_i^j}$ be the element of the canonical base of $\Z^\gamma$ associated to that vertex according to that order and $r_i^j$ the relation given by the multidegree of $v_i^j$.
We claim that the DCG of $CS^n(\underline{k}, \underline{h})$, $\Z^\gamma/<r_i^j>$, is generated by $[e_{v_2^1} - e_{v_1^1}]$.
Let $G:=<[e_{v_2^1} - e_{v_1^1}]>$. As $G$ is a subgroup of $\Z^\gamma/<r_i^j>$, it is a finite cyclic group. We shall prove that $G=\Z^\gamma/<r_i^j>$.
As $\Z^\gamma$ is generated by $\{e_{v_i^j} - e_{v_1^1}, i= 1,\ldots k_j, j=1,\ldots ,n\}$, the strategy will be to prove that every equivalence class $[e_{v_i^j}-e_{v_1^1}]\in G, i=1, \ldots,k_j, \, j=1,\ldots , n$.
Let us proceed by induction on $j$.

First we shall prove that all equivalence classes $[e_{v_i^1}-e_{v_1^1}]\in G, i=1, \ldots,k_1$.
To simplify the proof, we will consider that $h_1=k_1-1$ and proceed in 2 steps:
\begin{enumerate}
\item For $ i\le h_1=k_1-1$;
\item For $i=k_1$.
\end{enumerate}

\noindent (1) If $h_1=1$, there is nothing to prove.
If $h_1\geq 2$, we will again proceed by induction, this time on $i$.
For $i=2$, that $[e_{v_2^1}-e_{v_1^1}]\in G$ is just the hypothesis.
Now, for $i>2$, suppose that $[e_{v_s^1}-e_{v_1^1}]\in G$ for $1\leq s\leq i$.
 If $i=h_1$, it is done.
If $i<h_1$, then the vertex $v_i^1$ has degree $2$: it is adjacent to $v_{i-1}^1$ and to $v_{i+1}^1$. 
So, $r_i^1=e_{v_{i-1}^1} - e_{v_i^1} + e_{v_{i+1}^1}-e_{v_i^1}$, and, using the equality 
$$G\ni 2[e_{v_i^1}-e_{v_1^1}]=[e_{v_i^1}-e_{v_1^1}+r_i^1]=[e_{v_{i-1}^1}-e_{v_1^1}]+[e_{v_{i+1}^1}-e_{v_1^1}]$$
and the inductive hypothesis, we conclude that
$[e_{v_{i+1}^1}-e_{v_1^1}]\in G$.

\medskip

\noindent (2) If $k_1=2$, then we're done. If $k_1>2$, then the vertex $v_1^1$ has degree 2: it is attached to $v_2^1$ and also to $v_{k_1}^1$. So, $r_1^1=e_{v_2^1}-e_{v_1^1}+e_{v_{k_1}^1}-e_{v_1^1}$ and we get
$$\left[-(e_{v_2^1}-e_{v_1^1})+r_1^1\right]=[e_{v_{k_1}^1}-e_{v_1^1}],$$
which allows us to conclude that $[e_{v_{k_1}^1}-e_{v_1^1}]\in G$.
So, $[e_{v_i^1}-e_{v_1^1}]\in G, i=1, \ldots,k_1$.

\medskip

Now, admitting that $[e_{v_i^s}-e_{v_1^1}]\in G$ for $i=1,\ldots, k_s$ and $s=1,\ldots, j$, we shall prove that $[e_{v_i^{j+1}}-e_{v_1^1}]\in G$ for $i=1,\ldots, k_{j+1}$.
Again, the proof will be made in 2 steps:
\begin{enumerate}
\item For $1\le i\le h_{j+1}$
\item For $h_{j+1} < i\le k_{j+1}$.
\end{enumerate}

\noindent (1) If $h_{j+1}=1$, as $v_1^{j+1}= v_{h_j}^j$, then we're done.
If $h_{j+1}>1$, then we use induction on $i$, $2\le i\le h_{j+1}$.
So, first we shall prove that $[e_{v_2^{j+1}}-e_{v_1^1}]\in G$. As $h_{j+1}\ge 2$, $v_1^{j+1}$ is not adjacent to any vertex of the type $v_i^k$, for $k>j+1$ (except of course in the case $h_{j+1}=2$, when $v_{h_{j+1}}^{j+1}= v_1^{j+2}$). More, it has exactly two edges to vertices $v_i^{j+1}$:
$$v_1^{j+1}\longleftrightarrow v_2^{j+1} \mbox{ and }
v_{h_j}^j= v_1^{j+1}\longleftrightarrow v_{k_{j+1}}^{j+1}= v_{h_j+1}^j.$$ 
Observe also that
$$r_1^{j+1}=\sum_{v\,adj \, to \, v_1^{j+1}}(e_v - e_{v_1^{j+1}}).$$
So, if $n$ is the degree of the vertex $v_1^{j+1}$, we have:
\begin{equation*}
\begin{split}
G\ni n[e_{v_{h_j}^j}-e_{v_1^1}] & =  n[e_{v_1^{j+1}}-e_{v_1^1}]=[n e_{v_1^{j+1}}- n e_{v_1^1} + r_1^{j+1}] \\
& = \left[ne_{v_1^{j+1}}-ne_{v_1^1}+\sum_{v\,adj \, to \, v_1^{j+1}}(e_v - e_{v_1^{j+1}})\right] \\
& = \left[\sum_{v\,adj \, to \, v_1^{j+1}}(e_v - e_{v_1^1})\right] \\
& = \left[e_{v_2^{j+1}}-e_{v_1^1}\right]+\left[\sum_{\substack{v\,adj \, to \, v_1^{j+1} \\ v\neq v_2^{j+1}}} (e_v - e_{v_1^{j+1}})\right].
\end{split}
\end{equation*}
But we already know that, except $v_2^{j+1}$, $v_1^{j+1}$ is only adjacent to vertices of the type $v_i^s$, with $s\le j$. So, $$\left[\sum_{\substack{v\,adj \, to \, v_1^{j+1} \\ v\neq v_2^{j+1}}} (e_v - e_{v_1^{j+1}})\right]\in G \Rightarrow [e_{v_2^{j+1}}-e_{v_1^1}]\in G.$$
Now, suppose that $[e_{v_i^{j+1}}-e_{v_1^1}] \in G$ for $2\le s\le i$. If $i= h_{j+1}$, then we're done. If not, by the equality
\begin{equation}
G \ni 2[e_{v_i^{j+1}}-e_{v_1^1}]  = [2 e_{v_i^{j+1}}-2e_{v_1^1}+r_i^{j+1}] 
 = [e_{v_{i+1}^{j+1}}-e_{v_1^1}] + [e_{v_{i-1}^{j+1}}-e_{v_1^1}], 
\end{equation}
and by the inductive hypothesis, we conclude that $[e_{v_{i+1}^{j+1}}-e_{v_1^1}]\in G$.

\medskip

\noindent(2) The procedure is analogous: we should start from the vertex $v_{k_{j+1}}^{j+1}=v_{h_j+1}^{j}$ and advance in the opposite direction untill we reach $v_{h_{j+1}+1}^{j+1}$.
\end{proof}

Although Theorem \ref{ciclici} describes a whole family of graphs having cyclic DCG, they are not the only ones with this property. In fact, other examples can be obtained using Theorem \ref{dollaro}. 

\subsubsection*{A formula for the complexity of $CS^n(\underline k,\underline h)$}

For $n=1$, $CS^1(k,h)=\mathcal C_k$, so its complexity is $k$. 
For $n=2$, $CS^2((k_1,k_2),(h_1,h_2))$ is made of two cycles of order $k_1$ and $k_2$ attached in one edge. As it can be easily seen directly, or applying Proposition \ref{comb} to any edge except the common one, its cardinality is $k_1k_2-1$. 
For $n=3$, applying again formula \ref{fond} to any edge $l$ of the third cycle $\mathcal C_{k_3}$ (except the one in common with the second cycle), we get:
$$
\begin{array}{rl}
c(CS^3((k_1,k_2,k_3),(h_1,h_2,h_3)))&=c(CS^3((k_1,k_2,k_3),(h_1,h_2,h_3))\cdot l)+\\
&\quad +c(CS^3((k_1,k_2,k_3),(h_1,h_2,h_3))-l)\\
&=c(CS^3((k_1,k_2,k_3-1),(h_1,h_2,h_3)))+\\
&\quad +c(CS^2((k_1,k_2),(h_1,h_2))).
\end{array}
$$
Now, if $k_3-1\ge 2$, we can apply the same argument to $CS^3((k_1,k_2,k_3-1),(h_1,h_2,h_3))$ and we get
$$
\begin{array}{rl}
c(CS^3((k_1,k_2,k_3),(h_1,h_2,h_3)))&=c(CS^3((k_1,k_2,k_3-2),(h_1,h_2,h_3)))+\\
&\quad +2c(CS^2((k_1,k_2),(h_1,h_2))).
\end{array}
$$
By induction we obtain
$$
\begin{array}{rl}
c(CS^3((k_1,k_2,k_3),(h_1,h_2,h_3)))&=(k_3-2)c(CS^2((k_1,k_2),(h_1,h_2)))+\\
&\quad +c(CS^3((k_1,k_2,2),(h_1,h_2,h_3))).
\end{array}
$$
Observe that
$$
\begin{array}{rl}
c(CS^3((k_1,k_2,2),(h_1,h_2,h_3))\cdot v_1^3v_2^3)&=c(CS^2((k_1,k_2),(h_1,h_2))\cdot v_{h_{2}}^{2}v_{h_{2}+1}^{2})\\
&=c(CS^2((k_1,k_2),(h_1,h_2)))-\\
&\quad-c(CS^2((k_1,k_2),(h_1,h_2))-v_{h_{2}}^{2}v_{h_{2}+1}^{2})\\
&=c(CS^2((k_1,k_2),(h_1,h_2)))-\\
&\quad-c(CS^1((k_1),(h_1)));
\end{array}
$$ 
so, the last step gives
$$
\begin{array}{rl}
c(CS^3((k_1,k_2,k_3),(h_1,h_2,h_3)))&=k_3 c(CS^2((k_1,k_2),(h_1,h_2)))-c(CS^1((k_1),(h_1)))\\
&=k_1k_2k_3-k_1-k_3. 
\end{array}
$$
In general, arguing the same way, we obtain 
\begin{prop}
The complexity of the graphs $CS^n$ is given by the following recursive formula 
\begin{equation}\label{kh}
\begin{array}{rl}
c(CS^n(\underline k,\underline h))&=k_n c(CS^{n-1}(k_1,\ldots,k_{n-1}),(h_1,\ldots,h_{n-1}))\\
& -c(CS^{n-2}((k_1,\ldots,k_{n-2}),(h_1,\ldots,h_{n-2}))).
\end{array}
\end{equation}
\end{prop}
Observe that this formula implies in particular (by induction) that \emph{$c(CS^n(\underline k,\underline h))$ depends only of $\underline k$ and not of $\underline h$.}

\bigskip

We can make a slightly more explicit computation when $k_{1}=\dots =k_{n}=k$. 
In this case $c_{n}(k):= c(CS^n(\underline k,\underline h))$ is a polynomial in $k$. Let
$$P_{n}(k):=\sum _{i=0}^{[n/2]}(-1)^{i}a_{n-2i}^{i}k^{n-2i}$$
be a polynomial of degree $n$ in $k$ defined recursively as follows:
$$a^{0}_{l}:=1\mbox{, }\forall l\geq 0$$
$$a^{m}_{l}:=\sum _{k=0}^{l}a^{m-1}_{k}\mbox{, }m\geq 1.$$
We assert that $P_{n}=c_{n}$ as polynomials, for any $n\geq 1$.
We prove this by induction on $n$.
For $n=1$, it's clear.
For $n=2$, $P_{2}(k)=k^{2}-a_{0}^{1}$. Since $a_{0}^{1}=a_{0}^{0}=1$,  $P_{2}(k)=k^{2}-a_{0}^{1}=k^{2}-1=c_{2}(k)$.
Suppose now that  $n\geq 2$, and that $c_{j}=P_{j}$ for any $j<n$.
Then, using formula \ref{kh} and the definition of the polynomial, we obtain the following equalities
$$
\begin{array}{rcl}
c_{n}(k)&=&kc_{n-1}(k)-c_{n-2}(k)\\
&=&kP_{n-1}(k)-P_{n-2}(k)\\
&=&\left(\sum _{i=0}^{[(n-1)/2]}(-1)^{i}a_{n-2i-1}^{i}k^{n-2i-1}\right)-\left(\sum _{h=0}^{[n/2]-1}(-1)^{h}a_{n-2(h+1)}^{h}k^{n-2(h+1)}\right)\\
&=&\sum _{i=0}^{[n/2]}(-1)^{i}b_{n-2i}^{i}k^{n-2i}
\end{array}
$$
where
$b_{n}^{0}=a_{n-1}^{0}=1$
and
$b_{n-2i}^{i}=a_{n-2i-1}^{i}+a_{n-2i}^{i-1}
=(\sum _{k=0}^{n-2i-1}a^{i-1}_{k})+a_{n-2i}^{i-1}
=a_{n-2i}^{i}$.
So we're done.


\subsection{List of graphs for $\overline M_2$ and $\overline M_3$}
Recall that a stable curve $C$ over $k$ is a nodal curve of genus $g\geq 2$ such that if $E\subset C$ is a smooth rational component, then $|E\cap\overline{C\setminus E}|\geq 3.$
Clearly this combinatorial condition on stable curves implies that there are only \emph{finitely many} possible graphs for stable curves of a fixed genus.
Next, we list all the possible graphs for stable curves of genus 2 and 3, as well as their complexity and their DCG structure. We will use $\mathbb{Z}_n$ to denote the quotient group $\mathbb{Z}/n\mathbb{Z}$. The graphs are ordered by increasing the number of nodes.
In the graphs we will indicate the geometric genus of each irreducible component only if it is not zero.

\newpage

\begin{itemize}
\item \textbf{Genus 2}


\begin{equation} \label{genus2}
\begin{tabular}{||c|c|c|c|c||} 
\hline \hline 
Graph configuration & Nodes & Components & Complexity & DCG \\
\hline
\xymatrix@=.2pc{
{}\\
*{\bullet} \ar@{}^<{2}  
} 
& 
\xymatrix@=.001pc{
\\
0}
&
\xymatrix@=.001pc{
\\
1}
& 
\xymatrix@=.001pc{
\\
1}
&
\xymatrix@=.001pc{
\\
0}\\
\xymatrix@=.2pc{
{}\\
*{\bullet} \ar@{-}^<{1} @(ur,dr)\\{}
}
&
\xymatrix@=.001pc{
\\
1}
&
\xymatrix@=.001pc{
\\
1}
&
\xymatrix@=.001pc{
\\
1}
&
\xymatrix@=.001pc{
\\
0}\\

\xymatrix@R=.1pc{
*{\bullet} \ar@{-}^<{1}^>{1}[r] & *{\bullet}
\\{}
}
&
1
&
2
&
1
&
0\\

\xymatrix@R=.2pc{
*{\bullet} \ar@{-} @(ur,dr) \ar@{-} @(ul,dl)
\\{}}
&
2
&
1
&
1
&
0\\

\xymatrix@R=.2pc{
*{\bullet} \ar@{-} @(ul,dl) \ar@{-}^>{1}[r] & *{\bullet}
\\{}}
&
2
&
2
&
1
&
0\\

\xymatrix@R=.2pc{
*{\bullet} \ar@{-}@(ul,dl)[] \ar@{-}[r] & *{\bullet} \ar@{-} @(ur,dr)
\\{}}
&
3
&
2
&
1
&
0\\

\xymatrix@R=.2pc{
*{\bullet} \ar @{-} @/_/[r] \ar@{-}[r] & *{\bullet} \ar@{-} @/_/[l]
\\{}}
&
3
&
2
&
3
&
$\mathbb Z_3$\\
\hline \hline
\end{tabular}
\end{equation}

\bigskip

\bigskip

\item \textbf{Genus 3}

\begin{equation*}
\begin{tabular}{||c|c|c|c|c||}
\hline \hline
Graph configuration & Nodes & Components & Complexity & DCG\\
\hline
\xymatrix@=.2pc{
{}\\
*{\bullet} \ar@{}^<{3}  
} 
& 
\xymatrix@=.001pc{
\\
0}
&
\xymatrix@=.001pc{
\\
1}
& 
\xymatrix@=.001pc{
\\
1}
&
\xymatrix@=.001pc{
\\
0}\\
\xymatrix@=.001pc{
\\
*{2 \,\bullet} \ar@{-}@(ur,dr) & {}
}
&
\xymatrix@=.001pc{
\\
1}
&
\xymatrix@=.001pc{
\\
1}
&
\xymatrix@=.001pc{
\\
1}
&
\xymatrix@=.001pc{
\\
0}
\\
&&&&\\
\xymatrix{
*{2\, \bullet} \ar@{-}[r] & *{\bullet \, 1}
}
&
1&2&
1
&
0\\
&&&&\\
\xymatrix{
*{\bullet} \ar@{-}_<{1}@(ul,dl) \ar@{-}@(dr,ur)
}
&
2&1&
1
&
0
\\
&&&&\\

\xymatrix{
*{\bullet} \ar@{-}@(ul,dl) \ar@{-}[r]^<{1}^>{1}& *{\bullet}
}
&
2&2&
1
&
0\\
&&&&\\
\xymatrix{
*{\bullet} \ar@{-}@(ul,dl) \ar@{-}[r]^>{2}& *{\bullet}
}
&
2&2&
1
&
0
\\
&&&&\\
\xymatrix{
*{\bullet} \ar @{-} @/_/[r]_>{1}  & *{\bullet} \ar@{-} @/_/[l]_>{1}
}
&
2&2&
2
&
$\mathbb{Z}_2$\\
\xymatrix@R=1pc{
*{\bullet} \ar@{-}[r]^<{1} & *{\bullet}\ar@{-}[r]^<{1}^>{1} & *{\bullet}\\
{}
}
&
2&3&
1
&
0\\
\xymatrix@R=1pc{
*{\bullet} \ar@{-}@(ul,dl) \ar@{-}@(ur,dr) \ar@{-}@(ur,ul)\\
{}
}
&
3&1
&
1&
0\\
\xymatrix@R=1pc{
*{\bullet} \ar@{-}@(ur,ul) \ar@{-}@(dr,dl) \ar@{-}[r]^>{1}& *{\bullet}\\
{}
}
&
3&2&
1
&
0\\
\hline
\hline
\end{tabular}
\end{equation*}

\begin{equation*}
\begin{tabular}{||c|c|c|c|c||}
\hline \hline
Graph configuration & Nodes & Components & Complexity & DCG\\
\hline 
{}&&&&\\
\xymatrix{
*{\bullet} \ar@{-}@(ul,dl) \ar@{-}[r]_<{1} & *{\bullet} \ar@{-} @(ur,dr)
}
&
3&2&
1
&
0\\
&&&&\\
\xymatrix{
*{\bullet} \ar @{-}@(ul,dl) \ar @{-}@/_/[r] & *{\bullet} \ar @{-}_<{1}@/_/[l]
}
&
3&2&
2
&
$\mathbb Z_2$\\
&&&&\\
\xymatrix{
*{\bullet} \ar @{-} @/_/[r] \ar@{-}[r] & *{\bullet} \ar@{-}_>{1} @/_/[l]
}
&
3&2&
3
&
$\mathbb Z_3$\\
&&&&\\
\xymatrix@R=.5pc{
*{\bullet} \ar@{-}@(ul,dl) \ar@{-}[r] & *{\bullet}\ar@{-}[r]^<{1}^>{1} & *{\bullet}\\
{}
}
&
3&3
&
1
&
0\\
\xymatrix@R=.2pc{
*{\bullet} \ar@{-}[r]^<{1} & *{\bullet}\ar@{-}@(ul,ur) \ar@{-}[r]^>{1} & *{\bullet}\\
{}
}
&
3&3&
1&
0\\
\xymatrix{
 *{\bullet} \ar@{-}^<{1}[r]& *{\bullet} \ar @{-} @/_/[r] & *{\bullet} \ar@{-}_<{1}@/_/[l]
}
&
3&3&
2
&
$\mathbb Z_2$\\
&&&&\\
\xymatrix@R=1pc{
*{\bullet} \ar@{-}@(ur,ul) \ar@{-}@(dr,dl) \ar@{-}[r]& *{\bullet} \ar@{-}@(ur,dr)\\
{}
}
&
4&2&
1&
0\\
\xymatrix{
*{\bullet} \ar @{-}@(ul,dl) \ar @{-}@/_/[r] & *{\bullet} \ar @{-}@/_/[l] \ar @{-}@(ur,dr)
}
&
4&2&
2
&
$\mathbb Z_2$\\
&&&&\\
\xymatrix{
*{\bullet} \ar@{-}@(ul,dl) \ar @{-} @/_/[r] \ar@{-}[r] & *{\bullet} \ar@{-}@/_/[l]
}
&
4&2&
3
&
$\mathbb Z_3$\\
&&&&\\
\xymatrix{
*{\bullet}  \ar @{-} @/_.3pc/[r] \ar @{-}@/_.8pc/[r]/& *{\bullet} \ar@{-}@/_.3pc/[l] \ar@{-}@/_.8pc/[l]
}
&
4&2&
4
&
$\mathbb Z_4$\\
&&&&\\
\xymatrix@R=.5pc{
*{\bullet} \ar@{-}@(ul,dl) \ar@{-}[r] & *{\bullet}\ar@{-}[r]^<{1} & *{\bullet}\ar@{-}@(ur,dr)\\
{}
}
&
4&3
&
1
&
0\\
\xymatrix@R=.2pc{
*{\bullet} \ar@{-}[r] \ar@{-}@(ul,dl) & *{\bullet}\ar@{-}@(ul,ur) \ar@{-}[r]^>{1} & *{\bullet}\\
{}
}
&
4&3&
1&
0\\
\xymatrix@R=.1pc{
{} & *{\bullet} \ar@{-}@(ul,dl) \ar@{-}[r]& *{\bullet} \ar @{-} @/_/[r] & *{\bullet} \ar@{-}_<{1}@/_/[l]\\
{}
}
&4&3
&
2
&
$\mathbb Z_2$\\
\xymatrix@=1.2pc{
& *{\bullet} \ar@{-}[dl] \ar@{-}^<{1}[dr]& \\
*{\bullet} \ar@{-}[rr] \ar@{-}@/_/[rr]& & *{\bullet}
}
&
\xymatrix@=.3pc{
\\
4}
&
\xymatrix@=.3pc{
\\
3}
&
\xymatrix@=.3pc{
\\
5\\
{}
}
&
\xymatrix@=.3pc{
\\
{\mathbb Z_5}
}
\\
\xymatrix@R=.8pc{
{} & *{\bullet} \ar@{-}^<{1}[r]& *{\bullet} \ar@{-}[r] \ar @{-} @/_/[r] & *{\bullet} \ar@{-}@/_/[l]&
}
&
4&3&
3
&
$\mathbb Z_3$
\\
\xymatrix@=.5pc{
*{\bullet}\ar@{-}^<{1}[dr]&&&\\
&*{\bullet}\ar@{-}[rr]&&*{\bullet}\ar@{-}@(ur,dr)\\
*{\bullet} \ar@{-}_<{1}[ur]&&
}&
\xymatrix@=.3pc{
\\
4
\\
{}
\\
{}}
&\xymatrix@=.3pc{
\\
4}
&
\xymatrix@=.3pc{
\\
1}
&
\xymatrix@=.3pc{
\\
0}
\\
\xymatrix@R=.3pc{
*{\bullet} \ar@{-}[r] \ar@{-}@(ul,dl) & *{\bullet}\ar@{-}@(ul,ur) \ar@{-}[r] & *{\bullet}
\ar@{-}@(ur,dr)\\
{}
}
&
5&3&
1&
0\\
\xymatrix@R=.3pc{
{}&&\\
*{\bullet} \ar@{-}[r] \ar@{-}@(ul,dl) & *{\bullet} \ar@{-}@/^/[r] \ar@{-}@/_/[r] & *{\bullet}
\ar@{-}@(ur,dr)\\
{}
}
&
\xymatrix@R=.3pc{
{}\\
5}
&
\xymatrix@R=.3pc{
{}\\
3}
&
\xymatrix@R=.3pc{
{}\\
2}
&
\xymatrix@R=.3pc{
{}\\
\mathbb{Z}_2}\\
\hline \hline
\end{tabular}
\end{equation*}

\begin{equation*}
\begin{tabular}{||c|c|c|c|c||}
\hline \hline
Graph configuration & Nodes & Components & Complexity & DCG\\
\hline 
{}&&&&\\
\xymatrix@R=.2pc{
{}\\
{} & *{\bullet} \ar@{-}[r] \ar@{-}@(ul,dl)& *{\bullet} \ar@{-}[r] \ar @{-} @/_/[r] & *{\bullet} \ar@{-}@/_/[l]&\\&&&&
}
&
\xymatrix@R=.2pc{
{}\\
5}
&
\xymatrix@R=.2pc{
{}\\
3}
&
\xymatrix@R=.2pc{
{}\\
3}
&
\xymatrix@R=.2pc{
{}\\
\mathbb Z_3}
\\
\xymatrix@=1.2pc{
& *{\bullet} \ar@{-}[dl] \ar@{-}@/_/[dl] \ar@{-}@/^/[dr] \ar@{-}[dr]& \\
*{\bullet} \ar@{-}[rr] & & *{\bullet}\\{}
}
&
\xymatrix@=.3pc{
\\
5}
&
\xymatrix@=.3pc{
\\
3}
&
\xymatrix@=.3pc{
\\8
\\
{}
}
&
\xymatrix@=.3pc{
\\
{\mathbb Z_8}}\\
\xymatrix@=1.2pc{
& *{\bullet} \ar@{-}[dl] \ar@{-}[dr] \ar@{-}@(ul,ur)& \\
*{\bullet} \ar@{-}[rr] \ar@{-}@/_/[rr]& & *{\bullet}
}
&
\xymatrix@=.5pc{
\\5}
&
\xymatrix@=.5pc{
\\3
\\
{}
}
&
\xymatrix@=.5pc{
\\5
\\
{}
}
&
\xymatrix@=.5pc{
\\
{\mathbb Z_5}}
\\
\xymatrix@=1.2pc{
&*{\bullet} \ar@{-}^<{1}[d]&\\
& *{\bullet} \ar@{-}[dl] \ar@{-}[dr] & \\
*{\bullet} \ar@{-}[rr] \ar@{-}@/_/[rr]& & *{\bullet}
}
&\xymatrix@=.5pc{
{}\\
{}\\
5
}
&
\xymatrix@=.5pc{
{}\\
{}\\
4}
&
\xymatrix@=.5pc{
{}\\
{}\\
5
\\
{}\\{}
}
&
\xymatrix@=.5pc{
{}\\
{}\\
{\mathbb Z_5}}\\
\xymatrix@R=.2pc{
{} & *{\bullet} \ar@{-}[r] \ar@{-}@(ul,dl)& *{\bullet} \ar @{-} @/_/[r] & *{\bullet} \ar@{-}@/_/[l] \ar@{-}^>{1}[r]& *{\bullet}
\\&&&&
}
&5&4&
2
&
$\mathbb Z_2$
\\
\xymatrix@=.5pc{
*{\bullet}\ar@{-}^<{1}[ddd]&&&\\
&&*{\bullet}\ar@{-}@(ur,dr)&\\
&&&\\
*{\bullet}\ar@{-}[rrr]\ar@{-}[rruu]&& &*{\bullet}\ar@{-}@(ur,dr)\\
{}
}
&\xymatrix{
\\5}
&\xymatrix{
\\4}
&
\xymatrix{
\\1}
&
\xymatrix{
\\
0}
\\
\xymatrix@R=.2pc{
*{\bullet} \ar@{-}[r] \ar@{-}@(ul,dl)& *{\bullet} \ar @{-} @/_/[r] & *{\bullet} \ar@{-}@/_/[l] \ar@{-}[r]& *{\bullet} \ar@{-}@(ur,dr)
\\{}
}
&6&4&
2
&
$\mathbb Z_2$
\\
\xymatrix@=1pc{
*{\bullet} \ar@{-}[dd]\ar@{-}@/_.5pc/[dd]\ar@{-}[drr]&&&\\
&&*{\bullet}\ar@{-}[r]&*{\bullet} \ar@{-}@(ur,dr)\\
*{\bullet} \ar@{-}[urr]&&&
}
&
\xymatrix@=.6pc{
\\6}
&\xymatrix@=.6pc{
\\4}
&
\xymatrix@=.6pc{
{}\\
5
\\
{}\\{}
}
&
\xymatrix@=.6pc{
{}\\
{\mathbb Z_5}}\\
\xymatrix@=.5pc{
&*{\bullet}\ar@{-}[ddd] \ar@{-}@(ur,ul)&\\
*{\bullet} \ar@{-}@(ul,dl)&&*{\bullet}\ar@{-}@(ur,dr)&\\
&&&\\
&*{\bullet}\ar@{-}[ruu]\ar@{-}[luu] & \\
{}
}
&\xymatrix@=.5pc{
\\6}
&\xymatrix@=.5pc{
\\4}
&
\xymatrix@=.5pc{
\\1}
&
\xymatrix@=.5pc{
\\
0}
\\
\xymatrix{
*{\bullet} \ar@{-}[r] \ar@{-}[d] \ar@{-}@/_/[d] & *{\bullet} \ar@{-}[d] \ar@{-}@/^/[d]\\
*{\bullet} \ar@{-}[r] & *{\bullet}
}
&\xymatrix@=.5pc{
\\6}
&\xymatrix@=.5pc{
\\4}
&
\xymatrix@=.5pc{
\\12}
&
\xymatrix@=.5pc{
\\
\mathbb Z_{2}\times \mathbb Z_6
}
\\
\xymatrix@=1.2pc{
&*{\bullet}\ar@{-}[dd] \ar@{-}[dl] \ar@{-}[dr] &\\
*{\bullet} \ar@{-}[rr] |!{[r]}\hole & &*{\bullet}\\
&*{\bullet}\ar@{-}[ru]\ar@{-}[lu] & 
}
&\xymatrix@=.5pc{
\\6\\
{}}
&\xymatrix@=.5pc{
\\4}
&
\xymatrix@=.5pc{
\\16}
&
\xymatrix@=.5pc{
\\
\mathbb Z_4\times\mathbb Z_4}
\\
\hline \hline
\end{tabular}
\end{equation*}

\end{itemize}


\section{The behaviour of the DCG under standard geometrical operations on the curve}\label{sec3}

Applying standard geometrical operations to the nodes of a curve $C$,
such as the blow up or the normalisation, one gets a new curve $C'$. 
In this section we relate the DCG of $C'$
to the one of $C$. 
We will consider the following operations on the curve $C$: normalisation,
blow up and smoothing of a node $P$, and we will denote the new curves respectively
$B_PC$, $N_PC$ and $S_PC$.
For the geometric definitions of these operations, see for instance \cite{har}.

The operation of blow up of a node is defined in the context of algebraic geometry
using deformations of $C$, i.e. algebraic families of curves which have $C$ as a special fibre.
this operation depends on the choice of the family.
Hence, we will consider the following as the definition of the blow up of $C$
in a point $P$:

\begin{defi}
Let $C$ be a curve, $P$ a node of $C$.  
The blow up of $C$ in $P$, denoted $B_PC$, is the curve obtained
attaching a $\mathbb P^1$ to $N_PC$ by joining two distinct points to the preimages of $P$ in $N_PC$.
\end{defi}
\begin{rem}{\upshape 
Note that when we see $C$ as a special fibre of a one-parameter family with total space smooth, and we blow up the point $P$ corresponding to the node of $C$, the new fibre $C'$ in the blown up family does \emph{not} correspond to the modification described above, as the exceptional $\mathbb P^1$ has multiplicity $2$ in $C'$. 
Hence,  in order to treat this case, one needs to generalise the notion of DCG to curves with multiple components, as done for instance in \cite{BLR} and in \cite{Lorarit} (in this last paper  there is precisely the combinatorial description of the geometrical blow up of a smooth family). 
If, on the other hand,  we consider a one-parameter family with  a  rational singularity at $P$, blowing up $P$ we obtain as new fibre exactly the one we describe in the above definition. }
\end{rem}

The corresponding modification of the topological structure are reflected in the dual graph as follows: 
Let $P$ be a node of $C$ and call $l$ the corresponding edge in $\Gamma_C$.
\begin{itemize}
\item to take the {\em normalisation} $N_PC$ of $C$ in $P$ corresponds to {\em deleting the edge $l$} in $\Gamma_C$;
\item to {\em blow up $C$ in $P$}, denoted $B_PC$, corresponds to substituting $l$ with two edges $p,q$ and a new vertex $v$ as in figure 
\ref{blowup};

\medskip

\begin{figure}[!htp]
\centerline{
\xymatrix@R=1pc{
 & *{\bullet} \ar@{-}[rr]|-{ \,l\, } & & *{\bullet} & \ar@{:>}[rr] & &  &  
*{\bullet} \ar@{-}[r]|-{p} & *{\bullet} \ar@{-}[r]|-{q} & *{\bullet}\\
\ar@{.}[ur] & & & & \ar@{.}[ul] & &   \ar@{.}[ur] &&&& \ar@{.}[ul]
}
}
\caption{Blow up.}\label{blowup}
\end{figure}

\item to take the {\em smoothing} $S_PC$ of $C$ in $P$ corresponds to {\em contracting $l$} in $\Gamma_C$, i.e. to identify the vertices that contain it.
\end{itemize}

Our key tool will be formula (\ref{fond}) given in Proposition \ref{comb}.
A first geometric interpretation of this formula follows directly from the observations made above:
{\em if $C$ is a curve and $P\in C$ is a node which connects two different components of $C$, then 
$$c(C)=c(N_PC)+c(S_PC), $$
i.e. the complexity of $C$ is equal to the complexity of its normalisation at $P$
plus the complexity of its smoothing at $P$.}

\subsection{Blow up and normalisation}\label{blow}

The following result is another translation of the equality (\ref{fond}) in terms of 
blow up and normalisation in a node.

\begin{prop}
Let $C$ be a curve and $P\in C$  a node which connects two different components of $C$, then 
\begin{equation}\label{NeB}
c(B_PC)=c(C)+c(N_PC).
\end{equation}
\end{prop}
\begin{proof} Call $l$ the edge associated to $P$ in the graph of $C$. 
Let $p$, $q$ be the new edges that substitute $l$ in  $\Gamma_{B_PC}$.
Applying equality (\ref{fond}) to $\Gamma_{B_PC}$ with $e=q$ (or equivalently $e=p$)
we get
$$c(\Gamma_{B_PC})=c(\Gamma_{B_PC}-q)+c(\Gamma_{B_PC}\cdot q).$$ 
Observe that $\Gamma_{B_PC}\cdot p=\Gamma_C$; on the other hand $\Gamma_{B_PC}-p$ is $\Gamma_C-l$ 
with a tail made of an edge and a vertex attached in a vertex , so clearly these two graphs have the same complexity.
\end{proof}

What happens when we perform the blow ups several times in more than one node?
We give here a general formula which answers to this question. 
Suppose first that we blow up $k$ times one node $P$ which connects two different components of $C$.
The result is the curve  obtained attaching a chain of $k$ rational components to the preimages of $P$ in $N_PC$.
By induction on $k$ it is easy to prove the following formula
\begin{equation}\label{kNeB}
c(B_{kP}C)=c(C)+kc(N_PC).
\end{equation}

Let us call $\{P_1, P_2,\cdots,P_\delta\}$ the set of nodes of $C$. 
Suppose that none of them joins the same irreducible component.
Let $\underline k=(k_1,k_2,\cdots,k_\delta)$ be a $\delta$-uple of nonnegative integers.
We will call $B_{\underline k}C$ the curve obtained performing $k_i$ blow ups on the node $P_i$
(notice that this curve doesn't depend on the order in which the successive blow ups are made).
Notice that $B_{\underline k}C=B_{k_1P_1}B_{k_2P_2}...B_{k_\delta P_\delta}C$, and that
$B_{\underline k}(B_{\underline h}C)=B_{\underline{k+h}}C$.
If $T\subseteq \{1,2,...,\delta\}$ we call $N_TC$ the normalisation of $C$ in all the nodes $P_i$, $i\in T$.

\begin{teo}\label{perepe}
With the above notations, if $S=\{i\in \Z \mid k_i\not = 0\}$
$$
c(B_{\underline k}C)=\sum_{T\subseteq S}\left(\prod_{i\in T}k_i\right)c(N_TC).
$$
\end{teo}
\begin{proof} We proceed by induction on $n=\sharp S$. 
When $n=1$ we are reduced to formula (\ref{kNeB}).
Let $n>1$. We can suppose that $S=\{1,2,...,n\}$.
Call $\underline k'=\underline k-k_n e_n$.
By induction hypothesis
$$
c(B_{\underline k'}C)=\sum_{\begin{matrix}T\subseteq S\\n\notin T \end{matrix}}\left( \prod_{i\in T}k_{i}\right)c(N_TC).
$$
Applying formula (\ref{kNeB}) to $B_{\underline k'}C$ with $k=k_n$ and $P=P_n$
and substituting the above relation,  we get
$$
c(B_{\underline k}C)=c(B_{k_nP_n}B_{\underline k'}C)=c(B_{\underline k'}C)+k_nc(N_{P_n}B_{\underline k'}C)=
c(B_{\underline k'}C)+k_nc(B_{\underline k'}N_{P_n}C)=
$$
$$
=\sum_{\begin{matrix}T\subseteq S\\n\notin T \end{matrix}}\left( \prod_{i\in T}k_{i}\right)c(N_TC)+
k_n\sum_{\begin{matrix}T\subseteq S\\n\notin T\end{matrix}}\left( \prod_{i\in T}k_{i}\right)c(N_{T\cup\{n\}}C),
$$
which is our claim.
\end{proof}
Observe that we can allow the summand to run over all subsets of $\{1,2,...,\delta\}$, since the additional terms are zero. 
When $P$ is a node contained in only one irreducible component of $C$, the corresponding edge is a loop. To blow up  $k$   times $P$ means to substitute in the graph the loop with a $k$-cycle. 
So the DCG turns out to have a new factor $\Z/k\Z$.

 
\subsection{The blow up of vine curves}\label{vine}

Let $D_{N}$ be a nodal curve which is union of two smooth curves $A$ and $B$ intersecting in $N$ nodes.
We will call such a curve a {\em vine curve}.

Let $\underline{m}$ be a $N$-uple of positive integers $m_{1}$,\dots $m_{N}$.
Call $D_N(\underline m)$ the blow up of $D_N$ $m_i-1$ times in the $i$-th node:
$D_{N}(\underline{m})\colon =B_{\underline m-1}D_{N}$. 
We can suppose $m_{i}\geq m_{i+1}$ for any $i$.  
When $m_{1}=\dots =m_{k}=m$ and $m_{k+1}=\dots =m_{N}=1$, we will call the resulting curve $D_{N}(k^{m})$.

In what follows we analyse the order and the structure of the DCG of $D_N(\underline m)$.
This problem has been completely solved in \cite{BLR} (prop. 10 of section 9.6), using a criterion of Bourbaki to determine the diagonal form of the intersection matrix.
In our approach the computation of the order is a simple application of Theorem \ref{perepe}. 
For what concerns the computation of the structure, we explicitly compute the order of a set of generators of the DCG of $D_N(\underline m)$.
In some cases, this implies that the DCG is cyclic. 
The same computation has been performed by Lorenzini in \cite{Lorlap} (example 2.5 and successive claims), using a more general method developed in the same article.

\subsubsection*{Degree class group order} 
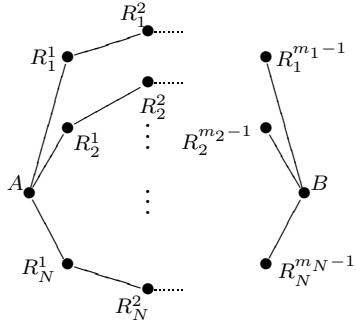
\begin{figure}[!htp]
\begin{equation*}
\xymatrix@R=.3pc@C=.7pc{
&&&*{\bullet} \ar@{.}[r] &&&&\\
&*{\bullet} \ar@{-}[urr]^>{R^2_1} &&&&&*{\bullet}\\
&&&*{\bullet} \ar@{.}[r]&&&\\
&*{\bullet} \ar@{-}[urr]_<{R^1_2}_>{R^2_2} &&{\vdots}&&&*{\bullet} \ar@{}[ul]^<{R^{m_2-1}_2} \\
*{\bullet} \ar@{-}[ur] \ar@{-}[uuur]^<{A}^>{R^1_1} \ar@{-}[ddr]_>{R^1_N} &&&{\vdots}&&&& *{\bullet} \ar@{-}[ul] \ar@{-}[uuul]_<{B}_>{R^{m_1-1}_1} \ar@{-}[ddl]^>{R^{m_N-1}_N} \\
\\
&*{\bullet} \ar@{-}[drr]_>{R^2_N} &&&&& *{\bullet} \\
&&&*{\bullet} \ar@{.}[r] &&&
}
\end{equation*}
\caption{$\Gamma_N(\underline{m})$} \label{blowvine}
\end{figure}

Let $\Gamma _{N}(\underline{m})$ be the dual graph of $D_{N}(\underline{m})$, which is composed of $N$ paths made of $m_{1},\dots m_{N}$ edges, such that every path links the vertex $A$ to the vertex $B$. (See Fig. \ref{blowvine}).

\begin{prop}\label{vineordine}
Let $c_{N}(\underline{m})$ be the complexity of $\Gamma _{N}(\underline{m})$.
$$c_{N}(\underline{m})=\overset{N}{\underset{k=1}{\sum }}\underset{i\neq k}{\prod }m_{i}$$
\end{prop}
\begin{proof}
Let $T\subseteq \{1,2,...,N\}$.
Observe that 
$N_TD_{N}=D_{N-\sharp T}$ 
and that
$$\sharp \{ T\subseteq T^{\prime }\subseteq \{1,2,...,N\} \mid \sharp T^{\prime }=N-1\}=N-\sharp T=c(N_TD_{N})\mbox{.}$$ 
By Proposition \ref{perepe} 
$$
c_{N}(\underline{m})=\sum_{T}\left(\prod_{i\in T}(m_i-1)\right)c(N_TC)\mbox{,}
$$
so
$$
\begin{array}{rl}
c_{N}(\underline{m})& =\sum_{T}\sum _{T\subseteq T^{\prime }\subseteq \{1,2,...,N\} }\left(\prod_{i\in T}(m_i-1)\right)=\\
 & =\overset{N}{\underset{k=1}{\sum }}\sum_{T\subseteq \{1,2,..,\hat{k},..,N\}}\left(\prod_{i\in T}(m_i-1)\right)=\\
& = \overset{N}{\underset{k=1}{\sum }}\underset{i\neq k}{\prod }(m_{i}-1+1).\\
\end{array}
$$
\end{proof}
\subsubsection*{Degree class group structure}

Let $\Delta _{N}(\underline{m})$ be the DCG of $D_{N}(\underline{m})$.
Let $R^{j}_{i}$ be the $j$-th component (from $A$ to $B$) of the chain associated to the $i$-th node. 
Also, we define $R^{0}_{i}\colon =A$ and $R^{m_{i}}_{i}\colon =B$, for every $i=1,\dots N$.
Suppose that $n_{0}$ is the maximum integer such that $m_{n_{0}}>1$.
Then the multidegrees of the components of $D_{N}(\underline{m})$ are the following:
$$
\begin{array}{rll}
\underline{c}_{A}=&-Ne_{A}+(N-n_{0})e_{B}+\overset{n_{0}}{\underset{i=1}{\sum }}e_{R^{1}_{i}}&=-Ne_{A}+\overset{N}{\underset{i=1}{\sum }}e_{R^{1}_{i}}\\ 
\underline{c}_{B}=&-Ne_{B}+(N-n_{0})e_{A}+\overset{n_{0}}{\underset{i=1}{\sum }}e_{R^{m_{i}-1}_{i}}&=-Ne_{B}+\overset{N}{\underset{i=1}{\sum }}e_{R^{m_{i}-1}_{i}}\\
\underline{c}_{R^{j}_{i}}=&-2e_{R^{j}_{i}}+e_{R^{j-1}_{i}}+e_{R^{j+1}_{i}}&\mbox{for }1\leq i\leq n_{0}\mbox{, }1\leq j\leq m_{i}-1\\
\end{array}
$$

\begin{prop}\label{vinegeneratori}

1) Let $t_{i}\colon =e_{R^{1}_{i}}-e_{A}$, $i=1,\dots ,N$. 
Then 
$$\Delta _{N}(\underline{m})=<t_{1},\dots ,t_{N}\mid \{ m_{i}t_{i}-m_{j}t_{j}\} _{i,j=1,\dots N}, 
\underset{i=1}{\overset{N}{\sum }}t_{i} >.$$
2) Let $M_{k}$ be the l.c.m. of the integers $m_{1},\dots \hat{m_{k}},\dots m_{N}$. 
Then the order of $[t_{k}]$ is 
$$\left(\underset{i=1}{\overset{N}{\sum }}(\underset{j\neq i}{\prod }m_{j})\right)(M_{k})/\left(\underset{j\neq k}{\prod }m_{j}\right)$$
\end{prop}
\begin{proof}
(1) According to Remark \ref{generatori}, $Z$ is generated by the elements $e_{R^{j}_{i}}-e_{A}$, $i=1,\dots N$, $j=0,\dots m_{i}$.
We will write $x\sim y$ instead of $[x]=[y]$.
We claim that 
$(e_{R^{j}_{i}}-e_{A})\sim jt_{i}$
for any $i=1,\dots n_{0}$ and $j=0,\dots m_{i}$.
In particular, $m_{i}t_{i}\sim e_{B}-e_{A}$ for any $i=1\dots n_{0}$, so we obtain the relations 
$m_{i}[t_{i}]=m_{k}[t_{k}].$

For fixed $i$, we will prove the claim by induction on $j$.
For $j=0$, it is clear.
Now, observe that
$$(e_{R^{j+1}_{i}}-e_{A})=2(e_{R^{j}_{i}}-e_{A})-(e_{R^{j-1}_{i}}-e_{A})+
\underline{c}_{R^{j}_{i}}\sim 2(e_{R^{j}_{i}}-e_{A})-(e_{R^{j-1}_{i}}-e_{A})$$ 
and suppose that the above claim is true for $j\leq l<m_{i}$, then 
$$e_{R^{l+1}_{i}}-e_{A}\sim 2(e_{R^{l}_{i}}-e_{A})-(e_{R^{l-1}_{i}}-e_{A})\sim 2lt_{i}-(l-1)t_{i}=(l+1)t_{i}.$$
Therefore,
$$\mbox{span}_{\Z }(t_{1},\dots , t_{n_{0}},e_{B}-e_{A})=Z/<\{ \underline{c}_{R^{j}_{i}}\} _{i\leq n_{0}, 0<j<m_{i}-1}>,$$
and $\underline{c}_{R^{m_{i}-1}_{i}}=(e_{B}-e_{A})-m_{i}t_{i}$, $i\leq n_{0}$.
Also, $\underline{c}_{A}=\overset{N}{\underset{i=1}{\sum }}t_{i}$ and $t_{i}=e_{B}-e_{A}$ for $i>n_{0}$, so
$$\mbox{span}_{\Z }(t_{1},\dots , t_{N})/<m_{i}t_{i}-m_{k}t_{k},\overset{N}{\underset{i=1}{\sum }}t_{i}>=
\Delta _{N}(\underline{m}).$$

\noindent (2) Let $d$ be a positive integer. 
Then $dt_{k}\sim 0$ if and only if $dt_{k}$ is a sum of multidegrees. 
Since $\underline{c}_{B}=\underset{I\neq B}{\sum }(-\underline{c}_{I})$, then $dt_{k}\sim 0$ if and only if 
$dt_{k}=\sum _{I\in V(D_{N}(\underline{m}))\setminus \{ B\} }a_{I}\underline{c}_{I}$ for some set of coefficients $a_{I}$, 
if and only if $d$ satisfies for some set of integers $a_{I}$ the following system

\begin{equation}\label{eq}
\left\{ \begin{array}{rclcr} 
-d&=&-Na_{A}+\overset{N}{\underset{i=1}{\sum }}a_{R^{1}_{j}}&\mbox{(projection on   }\Z e_{A}\mbox{    )}&\\

d&=&-2a_{R^{1}_{k}}+a_{A}+a_{R^{2}_{k}}&\mbox{       (projection on   } \Z e_{R^{1}_{k}}\mbox{ )}&\\
0&=&-2a_{R^{j}_{k}}+a_{R^{j-1}_{k}}+a_{R^{j+1}_{k}}&\mbox{       (projection on   }\Z e_{R^{j}_{k}}\mbox{ )} &j\neq 1\\
0&=&-2a_{R^{j}_{i}}+a_{R^{j-1}_{i}}+a_{R^{j+1}_{i}}&\mbox{       (projection on   }\Z e_{R^{j}_{i}}\mbox{ )} &i\neq k
\end{array}\right.
\end{equation}
Let $A^{j}_{i}\colon =a_{R^{j}_{i}}-a_{A}$. Observe that, for any $i$, $A^{0}_{i}=0$ and $A^{m_{i}}_{i}=a_{B}-a_{A}$.
Thus the system becomes
$$\left\{ \begin{array}{rclr} 
A^{1}_{k}+d&=&\underset{i\neq k}{\sum }(-A^{1}_{i})&\\
(A^{2}_{k}+d)&=&2(A^{1}_{k}+d)&\\
(A^{j+1}_{k}+d)&=&2(A^{j}_{k}+d)-(A^{j-1}_{k}+d)&j\neq 1\\
A^{j+1}_{i}&=&2A^{j}_{i}-A^{j-1}_{i}& i\neq k
\end{array}\right.$$
By the second and the subsequent equations follows that 
$$\begin{array}{rclr}
(A^{j}_{k}+d)&=&j(A^{1}_{k}+d)&\\
A^{j}_{i}&=&jA^{1}_{i} &i\neq k\\
\end{array}$$ 
for every $j$ (using induction on $j$!), so 
$$m_{k}(A^{1}_{k}+d)-d=A^{m_{k}}_{k}=(a_{B}-a_{A})=A^{m_{i}}_{i}=m_{i}A^{1}_{i}$$
for any $i\neq k$.
Then the previous system becomes 
$$\left\{ \begin{array}{rclr}
a_{B}-a_{A}&=&\underset{i\neq k}{\sum }(-m_{k})A^{1}_{i}-d&\\
a_{B}-a_{A}&=&m_{i}A^{1}_{i} &i\neq k
\end{array}\right.$$
The first equation by the product $\underset{j\neq k}{\prod }m_{j}$ becomes 
$$\begin{array}{rcl}
\underset{i\neq k}{\prod }m_{i}(a_{B}-a_{A})&=&\underset{i\neq k}{\sum }(-m_{k})(\underset{j\neq k}{\prod }m_{j})A^{1}_{i}-(\underset{j\neq k}{\prod }m_{j})d\\
\Lra \underset{i\neq k}{\prod }m_{i}(a_{B}-a_{A})&=&\underset{i\neq k}{\sum }(-\underset{j\neq i}{\prod }m_{j})(m_{i}A^{1}_{i})-(\underset{j\neq k}{\prod }m_{j})d\\
\Lra \underset{i\neq k}{\prod }m_{i}(a_{B}-a_{A})&=&-\left(\underset{i\neq k}{\sum }(\underset{j\neq i}{\prod }m_{j})\right)(a_{B}-a_{A})-(\underset{j\neq k}{\prod }m_{j})d\\
\Lra (\underset{j\neq k}{\prod }m_{j})d&=&\left(\underset{i=1}{\overset{N}{\sum }}(\underset{j\neq i}{\prod }m_{j})\right)(a_{A}-a_{B})
\end{array}$$
We now sum up the previous steps: 
if $dt_{k}\sim 0$ then there exists a set of integers $a_{I}$ such that $$d=\left(\underset{i=1}{\overset{N}{\sum }}(\underset{j\neq i}{\prod }m_{j})\right)(a_{A}-a_{B})/\left(\underset{j\neq k}{\prod }m_{j}\right)$$ 
and $a_{A}-a_{B}$ is a multiple of the $m_{i}$, $i\neq k$;
in particular a necessary condition for $d[t_{k}]=0$ is that $d$ is a multiple of $(\underset{i=1}{\overset{N}{\sum }}(\underset{j\neq i}{\prod }m_{j}))(M_{k})/(\underset{j\neq k}{\prod }m_{j})$. 

Vice versa, let $\overline{M}$ be any multiple of $M_{k}$ and let $a_{A}$ be any integer.
Then the following integers 
$$a_{B}\colon = a_{A}-\overline{M} \mbox{,  }$$ 
$$d\colon = \left(\underset{i=1}{\overset{N}{\sum }}(\underset{j\neq i}{\prod }m_{j})\right)(\overline{M})/\left(\underset{j\neq k}{\prod }m_{j}\right) \mbox{,  }$$ 
$$a_{R^{1}_{k}}\colon = -d+\underset{i\neq k}{\sum }(-a_{R^{1}_{i}}+a_{A}) \mbox{,  }$$ 
$$a_{R^{1}_{i}}\colon = \overline{M}/m_{i} \mbox{,  }i\neq k$$
and
$$a_{R^{j}_{k}}\colon = ja_{R^{1}_{k}}-(j-1)a_{A}+(j-1)d \mbox{,  }$$
$$a_{R^{j}_{i}}\colon = ja_{R^{1}_{i}}-(j-1)a_{A} \mbox{,  }i\neq k$$ 
for $j>1$, satisfy  system (\ref{eq}) at the beginning of the proof; 
in particular $dt_{k}\sim 0$ for $d\colon = (\underset{i=1}{\overset{N}{\sum }}(\underset{j\neq i}{\prod }m_{j}))(\overline{M})/(\underset{j\neq k}{\prod }m_{j})$, for every $\overline{M}$ multiple of $M_{k}$.
In conclusion, we have shown that $d[t_{k}]=0$ if and only if $d$ is a multiple of  $(\underset{i=1}{\overset{N}{\sum }}(\underset{j\neq i}{\prod }m_{j}))(M_{k})/(\underset{j\neq k}{\prod }m_{j})$. 
\end{proof}
\begin{cor}\label{coro}
 $\Delta _{N}(\underline{m})$ is generated by $[t_{k}]$ if and only if $(m_{i},m_{l})=1$ for all $i,l\neq k$.
\end{cor}
\begin{proof} Since the order of $\Delta _{N}(\underline{m})$ is $\underset{i=1}{\overset{N}{\sum }}(\underset{j\neq i}{\prod }m_{j})$ by Proposition \ref{vineordine} (and preliminaries) and the order of $[t_{k}]$ is $(\underset{i=1}{\overset{N}{\sum }}(\underset{j\neq i}{\prod }m_{j}))(M_{k})/(\underset{j\neq k}{\prod }m_{j})$ by Proposition \ref{vinegeneratori}, then $[t_{k}]$ generates $\Delta _{N}(\underline{m})$ if and only if $M_{k}=(\underset{j\neq k}{\prod }m_{j})$. 
\end{proof}

\begin{rem}
\upshape{
Observe that in general there exists no subset of $\{ t_{i}\} _{i=1,\dots N}$ such that the classes of its elements generate the cyclic factors of the DCG. 
Indeed, this is the case if and only if every 
non-trivial linear combination of such classes isn't zero, whereas in the first part of the proof of Proposition 
\ref{vinegeneratori} we obtained $m_{i}[t_{i}]+(-m_{j})[t_{j}]=0$, and $m_{i}$ is never zero.}
\end{rem}

\begin{rem}\label{presentation} 
\upshape{Observe that the presentation given in the previous proposition is equivalent to the following one:    
$$\Delta _{N}(\underline{m})=<t_{1},\dots ,t_{N-1}\mid \overset{N}{\underset{j=1}{\sum }}(m_{j}+\delta _{ij}m_{N})t_{j}>$$
We can rewrite this presentation as an exact sequence:
$$0\longrightarrow \Z ^{N-1}\stackrel{\Sigma}{\longrightarrow}\Z ^{N-1} \longrightarrow\Delta _{N}(\underline{m})\longrightarrow 0$$ 
where the endomorphism $\Sigma : \Z ^{N-1}\raw \Z ^{N-1}$ is represented, with respect to the canonical base,  by the matrix with entries 
$$a_{ij}=m_{N}+\delta _{ij}m_{i}.$$
The problem of the decomposition of the DCG in cyclic factors corresponds to  the problem of the decomposition of the coker of $\Sigma$; as it is well known, the latter problem is equivalent to the diagonalization of any matrix associated to $\Sigma$.}
\end{rem}

\begin{rem}\label{fattorecomune}
\upshape{Observe that if $m_{1},\dots m_{N}$ have a common factor $d$ and $\overline{m_{i}}\colon = m_{i}/d$, then $\Sigma$ is the 
composition of the multiplication by $d$, by the map $\overline{\Sigma }$ which is represented with respect to the canonical 
base by the matrix whose entries are
$$a_{ij}=\overline{m_{N}}+\delta _{ij}\overline{m_{i}}\mbox{ .}$$ 
So, the entries of the diagonal form of $\Sigma$ are $d$ times the entries of the diagonal form of $\overline{\Sigma }$, i.e. 
{\it if the coker of} $\overline{\Sigma }$ {\it is }$\underset{i=1}{\overset{s}{\bigoplus }}\Z /k_{i}\Z $
{\it , then the DCG is }
$$\underset{i=1}{\overset{s}{\bigoplus }}\Z /dk_{i}\Z \oplus (\Z /d\Z )^{N-s}$$}
\end{rem}
\noindent Thus we can suppose that $m_{1},\dots m_{N}$ have no common factor.
Using our results, we can compute the structure of $\Delta_N(\underline m)$ 
in the following cases:

\begin{prop}\label{muguali}
$$
\begin{array}{rcl}
\Delta _{N}(1^{m})&\cong&\Z /(1+m(N-1))\Z\mbox{,}\\
\Delta _{N}(k^{m})&\cong&\Z /m(k+m(N-k))\Z \oplus (\Z /m\Z )^{k-2}\mbox{, for }k\geq 2\end{array}
$$
\end{prop}
\begin{proof}
The first equality follows from Corollary \ref{coro}.
As for the second one, observe that the relations are generated by $mt_{2}-mt_{1},\dots ,mt_{k}-mt_{1}, t_{k+1}- mt_{1},\dots ,t_{N}- mt_{1}$ and $\underset{i=1}{\overset{N}{\sum }}t_{i}$. 
Thus we can forget the generators $t_{k+1},\dots ,t_{N}$, change the generators replacing $t_{i}=(t_{i}-t_{1})+t_{1}$, 
and obtain as relations the following ones:
$$m(t_{2}-t_{1}),\dots ,m(t_{k}-t_{1}),$$
and
$$\underset{i=2}{\overset{k}{\sum }}(t_{i}-t_{1})+(k+m(N-k))t_{1}.$$
Observe that $t_{k}-t_{1}$ belongs to the subgroup of the DCG which is generated by $t_{1}, t_{2}-t_{1},\dots ,t_{k-1}-t_{1}$,
 so we can delete it from the list of generators using the identity   
$$m(t_{k}-t_{1})=m(k+m(N-k))t_{1}+\underset{i=2}{\overset{k-1}{\sum }}m(t_{i}-t_{1})$$
and deleting the latter relation.
Therefore
$$\Delta _{N}(\underline{m})=<t_{1},t_{2}-t_{1},\dots ,t_{k-1}-t_{1}\mid m(t_{2}-t_{1}),\dots m(t_{k-1}-t_{1}), 
m(k+m(N-k))t_{1}>,$$   
and we are done.
\end{proof}

\subsubsection*{Structure of the DCG of the dollar sign curve}
The vine curve with three nodes, $\Delta _{3}(\underline{m})$, is usually called {\em dollar sign curve}
(cf. \cite{OS}, section 9; the reason is that the picture of the curve itself resembles to the dollar symbol).
Here, by means of an ad hoc algebraic argument, we show that the DCG of any iterated blow up of the dollar curve, $\Delta _{3}(\underline{m})$, is \av almost anytime'' a cyclic group; indeed, from the result below and Remark \ref{fattorecomune}, it follows that 
$$\Delta _{3}(\underline{m})\cong \mathbb Z/d\mathbb Z\otimes \mathbb Z/k\mathbb Z,$$
where $d=\gcd (\underline{m})$, and $k= \frac{c_3(\underline m)}{d}$ (where $c_3(\underline m)$ has been defined in Proposition \ref{vineordine}). 

\begin{prop}\label{dollaro}
The DCG $\Delta _{3}(\underline{m})$ of any (iterated) blow up of the dollar sign curve is a cyclic group whenever $\gcd (\underline{m})=1$.
\end{prop}
By Remark \ref{presentation},
the proposition is a special case of the following result:
\begin{lem}
Let $\Sigma$ be a endomorphism of $\Z ^{2}$ induced by a matrix $M$ of entries $m_{ij}$. 
Suppose that 
$$m_{11}>0\mbox{, }m_{12}\geq 0\mbox{ ( or }m_{21}\geq 0\mbox{ )}$$ 
and that
$$
\begin{array}{rcl}
\gcd (m_{21},\gcd (m_{11},m_{12}))&=&\gcd (m_{11},m_{12})\mbox{,}\\
\gcd (m_{22},\gcd (m_{11},m_{12}))&=&1\\
\end{array}
$$
( or, respectively
$$
\begin{array}{rcll}
\gcd (m_{12},\gcd (m_{11},m_{21}))&=&\gcd (m_{11},m_{21})&\mbox{,}\\
\gcd (m_{22},\gcd (m_{11},m_{21}))&=&1)&\\
\end{array}
$$

Then the coker of $\Sigma$ is a cyclic group.
\end{lem}
\begin{proof}
We will prove the thesis by induction on $m_{11}$ (which is a natural number by assumption).
For $m_{11}=1$,
$$ 
\begin{array}{rcl}
\left( \begin{array}{rl}
1&0\\
-m_{21}&1\\
\end{array}\right)
\left( \begin{array}{rl}
m_{11}&m_{12}\\
m_{21}&m_{22}\\
\end{array}\right)
\left( \begin{array}{rl}
1&-m_{12}\\
0&1\\
\end{array}\right)
&=& 
\left( \begin{array}{rl}
1&0\\
0&det(\Sigma )\\
\end{array}\right)
\end{array}
$$
Hence, the coker of $\Sigma$ is $\Z$ when $det(\Sigma )=0$, and it is $\Z /det(\Sigma )$ otherwise.  

For $m_{11}=m>1$, if $d\colon =\gcd (m_{11},m_{12})=m$ (and so $m_{21}$ is divisible by $m$).
Set $a\colon = m_{12}/m$, $b\colon = m_{21}/m$, then 
$$ 
\begin{array}{rcl}
\left( \begin{array}{rl}
1&0\\
-b&1\\
\end{array}\right)
\left( \begin{array}{rl}
m_{11}&m_{12}\\
m_{21}&m_{22}\\
\end{array}\right)
\left( \begin{array}{rl}
1&-a\\
0&1\\
\end{array}\right)
&=& 
\left( \begin{array}{rl}
m&0\\
0&m_{22}-am_{21}\\
\end{array}\right)
\end{array},
$$
so the coker of $\Sigma$ is decomposed as $\Z /m\Z \oplus \Z /(m_{22}-am_{21})\Z$. Since 
$$\gcd (m,m_{22}-am_{21})=\gcd (d,m_{22})=1,$$ 
we have 
$$\Z /m\Z \oplus \Z /(m_{22}-am_{21})\Z = \Z /det(\Sigma )\Z$$ 
and therefore the coker of $\Sigma$ is a cyclic group. 

On the other hand, suppose that $d<m$. We can write $d$ as a linear combination of $m_{11}$ and $m_{12}$:
$$d=x_{1}m_{11}+x_{2}m_{12}.$$ 
Let $m_{11}=a_{1}d$, $m_{12}=a_{2}d$. The matrix
$$
\begin{array}{rcl} 
A&\colon = &
\left( \begin{array}{rl}
x_{1}&-a_{2}\\
x_{2}&a_{1}\\
\end{array}\right)
\end{array} 
$$
is an integer invertible matrix, because 
$$\det A=x_{1}a_{1}+x_{2}a_{2}=(x_{1}m_{11}+x_{2}m_{12})/d=1.$$
In particular $a_{1}$ is prime with $x_{2}$.
The product $AM$ gives us a new matrix $M'$ associated to $\Sigma$, with the following entries
$$
\begin{array}{rcl}
m^{\prime }_{11}&\colon = &x_{1}m_{11}+x_{2}m_{12}=d>0,\\
m^{\prime }_{12}&\colon = &-a_{2}m_{11}+a_{1}m_{12}=-a_{2}a_{1}d+a_{1}a_{2}d=0,\\
m^{\prime }_{21}&\colon = &x_{1}m_{21}+x_{2}m_{22},\\
m^{\prime }_{22}&\colon = &-a_{2}m_{21}+a_{1}m_{22}.\\
\end{array}
$$
Note that
$$
\gcd (m^{\prime }_{12}, \gcd (m^{\prime }_{11},m^{\prime }_{21}))
= \gcd (0,\gcd (m^{\prime }_{11},m^{\prime }_{21}))
= \gcd (m^{\prime }_{11},m^{\prime }_{21}),
$$
and that
$$
\begin{array}{rcl}
\gcd (m^{\prime }_{22}, \gcd (m^{\prime }_{1,1},m^{\prime }_{21}))
&=& \gcd (-a_{2}m_{21}+a_{1}m_{22},\gcd (d,x_{1}m_{21}+x_{2}m_{22}))\\
&=& \gcd (-a_{2}m_{21}+a_{1}m_{22},\gcd (d,x_{2}))\\
&=& \gcd (\gcd (-a_{2}m_{21}+a_{1}m_{22},d),x_{2}))\\
&=& \gcd (\gcd (a_{1},d),x_{2}))\\
&=& \gcd (\gcd (d,a_{1}),x_{2}))\\
&=& \gcd (d,\gcd (a_{1},x_{2}))\\
&=& \gcd (d,1)=1
\end{array}
$$
Hence, the new matrix satisfies the hypothesis of the proposition, and $m^{\prime }_{11}=d<m$; and we can conclude by inductive hypothesis that the coker of $\Sigma$ is cyclic.
\end{proof}


\bibliographystyle{amsalpha}
  \nocite{*}
  \bibliography{artigo_prag2}

\providecommand{\bysame}{\leavevmode\hbox to3em{\hrulefill}\thinspace}
\providecommand{\MR}{\relax\ifhmode\unskip\space\fi MR }
\providecommand{\MRhref}[2]{%
  \href{http://www.ams.org/mathscinet-getitem?mr=#1}{#2}
}
\providecommand{\href}[2]{#2}
\begin{thebibliography}{BdlHN97}

\bibitem[Art91]{ar}
M.~Artin, \emph{Algebra}, Prentice Hall, 1991.

\bibitem[BdlHN97]{tatiana}
R.~Bacher, P.~de~la Harpe, and T.~Nagnibeda, \emph{The lattice of integral
  flows and the lattice of integral cuts on a finite graph.}, Bull. Soc. Math.
  France \textbf{125} (1997), no.~2, 167--198.

\bibitem[Ber70]{berge}
C.~Berge, \emph{Graphs and hypergraphs}, North-Holland, 1970.

\bibitem[Big74]{biggs}
N.~Biggs, \emph{Algebraic {G}raph {T}heory}, Cambridge University Press, 1974.

\bibitem[Big99]{biggschip}
\bysame, \emph{Chip-{F}iring and the {C}ritical {G}roup of a {G}raph}, Journal
  of Algebraic Combinatorics \textbf{9} (1999), 25--45.

\bibitem[BLR90]{BLR}
S.~Bosch, W.~L\"uktebohmert, and M.~Raynaud, \emph{N\'eron models}, Ergebnisse
  der Mathematik, no.~21, Springer-Verlag, 1990.

\bibitem[Cap]{prag}
L.~Caporaso, \emph{Introduction to moduli of curves}, Notes of the Summer
  School PRAGMATIC 2004.

\bibitem[Cap94]{cap}
\bysame, \emph{A compactification of the universal {P}icard variety over the
  moduli space of stable curves}, Journal of the American Mathematical Society
  \textbf{7} (1994), no.~3, 589--660.

\bibitem[Cap05]{capneron}
\bysame, \emph{N\'eron models over moduli of stable curves}, Preprint
  math.AG/0502171 (2005).

\bibitem[Cha82]{chaiken}
S.~Chaiken, \emph{A combinatoric proof of all minors matrix tree theorem}, SIAM
  Journal Algebraic Discrete Methods \textbf{3} (1982), 319.

\bibitem[Har77]{har}
R.~Hartshorne, \emph{Algebraic {G}eometry}, G.T.M., no.~52, Springer-Verlag,
  1977.

\bibitem[HM98]{harris}
J.~Harris and I.~Morrison, \emph{Moduli of {C}urves}, Springer-Verlag, 1998.

\bibitem[Lor89]{Lorarit}
D.~Lorenzini, \emph{Arithmetical graphs}, Math. Ann. \textbf{285} (1989),
  no.~3, 481--501.

\bibitem[Lor90a]{lorgraphs}
\bysame, \emph{Dual graphs of degenerating curves}, Mat. Annalen \textbf{287}
  (1990), 135--150.

\bibitem[Lor90b]{lorjac}
\bysame, \emph{Groups of components of {N}\'eron models of {J}acobians},
  Compositio Matematica \textbf{73} (1990), 145--160.

\bibitem[Lor91]{Lorgroup}
\bysame, \emph{A finite group attached to the {L}aplacian of a graph}, Discrete
  Math. \textbf{91} (1991), no.~3, 277--282.

\bibitem[Lor93]{lorner}
\bysame, \emph{On the group of components of a {N}\'eron model}, J. Reine
  Angew. Math. \textbf{445} (1993), 109--160.

\bibitem[Lor00]{Lorlap}
\bysame, \emph{Arithmetical properties of {L}aplacians of graphs}, Linear and
  Multilinear Algebra \textbf{47} (2000), no.~4, 281--306.

\bibitem[OS79]{OS}
T.~Oda and C.~Seshadri, \emph{Compactifications of the generalized {J}acobian
  variety}, Trans. A.M.S. \textbf{253} (1979).

\bibitem[Ray70]{ray}
M.~Raynaud, \emph{Specialisation du foncteur de {P}icard}, Inst. Hautes Etudes
  Sci. Publ. Math. \textbf{28} (1970), 27--76.

\bibitem[Roy01]{godsil-royle}
C.~Godsil;~G. Royle, \emph{Algebraic graphy theory}, Graduate Text in
  Mathematics, no. 207, Springer-Verlag, 2001.

\bibitem[Wes96]{west}
D.~West, \emph{Introduction to graph theory}, Prentice Hall, 1996.

\end{thebibliography}

\bigskip

\noindent Simone Busonero, Dipartimento di Matematica, Universit\`a di Roma la Sapienza. \\ 
E-mail: \textsl {busonero@mat.uniroma1.it}.
\bigskip

\noindent Margarida Melo, Departamento de Matem\'atica, Universidade de Coimbra - Dipartimento di Matematica, Unversit\`a di Roma Tre.\\
E-mail: \textsl {mmelo@mat.uc.pt, melo@mat.uniroma3.it}.

\bigskip
\noindent Lidia Stoppino, Dipartimento di Matematica, Unversit\`a di Roma Tre.\\ E-mail: \textsl {lidia.stoppino@unipv.it}.

\end{document}